\documentclass[11pt]{amsart}

\usepackage{amsmath,amssymb,amsthm}
\usepackage[margin=1.15in]{geometry}
\usepackage{enumerate}

\newtheorem{theorem}{Theorem}[section]
\newtheorem{mainthm}{Theorem}

\newtheorem{proposition}[theorem]{Proposition}
\newtheorem{lemma}[theorem]{Lemma}
\newtheorem{corollary}[theorem]{Corollary}
\theoremstyle{definition}
\newtheorem{definition}[theorem]{Definition}
\newtheorem{problem}[theorem]{Problem}
\theoremstyle{remark}
\newtheorem{remark}[theorem]{Remark}

\numberwithin{equation}{section}

\begin{document}

\title[Uniform estimates for the Erd\H{o}s--Moser equation]{Uniform estimates for the Erd\H{o}s--Moser equation\\ and its analogue in arithmetic progressions}

\author{Elham Mehdi-Nezhad}
\address{Department of Mathematics and Applied Mathematics, University of the Western Cape, Private Bag X17, Bellville 7535, Cape Town, South Africa}
\email{emehdinezhad@uwc.ac.za}

\author{Amir M. Rahimi}
\address{School of Mathematics, Institute for Research in Fundamental Sciences (IPM), P.O. Box 19395-5746, Tehran, Iran}
\email{amrahimi@ipm.ir}


\begin{abstract}
We study the Erd\H{o}s--Moser equation in arithmetic progressions, $\sum i^k=m^k$, where $i$ runs over the positive integers $i<m$ with $i\equiv b\pmod q$. We prove that a solution has at most $\lfloor k/\log2+\frac32\rfloor$ terms for all $k<10^{55}$, and at most $k/\log2+\frac52$ terms in general; this linear bound is best possible for real solutions and replaces the previously known exponential bound $(k+1)2^{2k+1}$. For $k=3$ we obtain a complete classification: solutions with exactly $n$ terms exist if and only if $n\in\{3,4,5\}$, and in each case there are infinitely many primitive solutions. For even exponents there are no solutions with three terms, and there are no solutions at all for $k=4$, for every $k=2^s\ge4$, and for every even $k$ with $4\le k\le24$. The key analytic input is a uniform estimate for the classical equation $1^k+\dots+(m-1)^k=m^k$: the real exponent $\kappa_m$ for which $1^{\kappa}+\dots+(m-1)^{\kappa}=m^{\kappa}$ satisfies $\kappa_m=\log2\,(m-\frac32-\gamma_m/m)$ with $0<\gamma_m<1/(8\log2)$ for every $m\ge3$. The proofs rest on a method based on centres of mass, Steiner's parallel-axis theorem and what we call the half-mass principle. The same framework identifies the constants of a general class of balance equations through the moments of a mass distribution: $\frac32$ is a Steiner quotient, and $\frac{25}{12}-3\log2=\frac1{12}-\frac12f''(\log2)$ with $f(\lambda)=\lambda/(e^\lambda-1)$, whose Bernoulli expansion explains why this constant is so small.
\end{abstract}

\maketitle

\section{Introduction}

\subsection{The equations}
Erd\H{o}s conjectured that the equation
\begin{equation}\label{eq:EM}
1^k+2^k+\dots+(m-1)^k=m^k
\end{equation}
has no solution in positive integers other than $1+2=3$. Moser \cite{Moser} showed that any further solution has $m>10^{10^6}$; Gallot, Moree and Zudilin \cite{GMZ} raised this bound to $m>2.7\cdot10^{10^9}$, using the asymptotic relation $k/m=\log2\,(1-\frac3{2m}-C_m/m^2)$ with $0<C_m<0.004$ for $m>10^9$ together with a large computation of the continued fraction of $\log2$; see the survey of Moree \cite{Moree} for the history. Recently Dong, Nguyen and Zaharescu \cite{DNZ} introduced the analogue in arithmetic progressions,
\begin{equation}\label{eq:APintro}
\sum_{\substack{1\le i\le m-1\\ i\equiv b\ (\mathrm{mod}\ q)}}i^k=m^k ,
\end{equation}
which, unlike \eqref{eq:EM}, has solutions: for instance $5^2+12^2=13^2$ with $q=7$. They described the solutions for $k=2$, bounded the number of terms by $(k+1)2^{2k+1}$ in general and by $10$ for $k=3$, and found solutions with three terms for $k=3$. This leaves three questions: how many terms a solution can have, which numbers of terms actually occur for a given exponent, and for which exponents solutions exist at all.

In this paper we replace the exponential bound for the number of terms in \eqref{eq:APintro} by a sharp linear bound, we classify completely the solutions for $k=3$, and we show that there are no solutions for $k=4$, and more generally for every $k=2^s\ge4$. The tool is the centre of mass; the same considerations explain the constants $\frac32$ and $\frac{25}{12}-3\log2$ of the Erd\H{o}s--Moser equation.

\subsection{Main results}
Our results combine a uniform estimate for the classical equation with applications to its analogue in arithmetic progressions. Let $\kappa_m$ be the real number for which $1^{\kappa}+\dots+(m-1)^{\kappa}=m^{\kappa}$; thus the solutions of \eqref{eq:EM} are the integers $m$ for which $\kappa_m$ is an integer, and we define $\gamma_m$ by $\kappa_m=\log2\,(m-\frac32-\gamma_m/m)$. Let $L=\log2$ and $\gamma_0=\frac{25}{12}-3L$.

\begin{mainthm}[Theorems~\ref{thm:main} and~\ref{thm:allm}]\label{thmA}
For every $m\ge3$ we have $0<\gamma_m<1/(8\log2)$, and for $m\ge10^5$ we have $-87/m<\gamma_m-\gamma_0<116/m$.
\end{mainthm}

The relation $\gamma_m\to\gamma_0$ is due to Gallot, Moree and Zudilin. The new point in Theorem~\ref{thmA} is uniformity: the bounds hold for every $m\ge3$ and use no arithmetic. This uniformity is what makes Theorem~\ref{thmA} applicable to equations with small solutions, such as \eqref{eq:APintro}; for \eqref{eq:EM} itself, whose non-trivial solutions are enormous, the asymptotic form of \cite{GMZ} suffices.

\begin{mainthm}[Proposition~\ref{prop:Gamma}, Theorem~\ref{thm:gen}, Corollary~\ref{cor:Rfamily}]\label{thmB}
Let $\Gamma$ be the system of masses $2^{-1},2^{-2},2^{-3},\dots$ placed at unit spacing to the left of a point $P$, let $d$ and $\rho$ be the distance of its centre from $P$ and its radius of gyration, and let $Q_1$ and $Q_2$ be its static moment and its moment of inertia about $P$. Then
\[
\frac32=\frac{\rho^2+d^2}{2d},\qquad \frac{25}{12}-3\log2=\frac1{12}+Q_1-\frac{\log2}{2}\,Q_2 .
\]
More generally, for the equation $1^k+\dots+(m-1)^k=Rm^k$ with $R>0$ and $\lambda=\log(1+1/R)$, the first two constants are $\frac12\coth\frac\lambda2$ and
\[
\frac1{12}-\frac12f''(\lambda)=\frac{\lambda^2}{120}-\frac{\lambda^4}{2016}+\frac{\lambda^6}{43200}-\cdots,\qquad f(\lambda)=\frac{\lambda}{e^\lambda-1}.
\]
\end{mainthm}

The coefficients of the last expansion are Bernoulli numbers, and it explains the near-cancellation $\frac{25}{12}\approx3\log2$: one has $(\log2)^2/120-(\log2)^4/2016<\gamma_0<(\log2)^2/120$. Theorem~\ref{thm:gen} gives the same description for a class of balance equations which includes \eqref{eq:APintro} for a fixed progression: the first constant is a Steiner quotient, and the second is given by the first three cumulants of a distribution of moments. This supplies the term of the expansion left aside in \cite[Remark~1.3]{DNZ}.

For the progression equation we normalise $1\le b\le q$ and let $n$ be the number of terms.

\begin{mainthm}[Lemma~\ref{lem:APreal}, Theorem~\ref{thm:APterms}]\label{thmC}
Let $k\ge2$. The equation \eqref{eq:APintro} with $n$ terms has real solutions $b/q$, $m/q$ in the admissible range if and only if $1^k+\dots+(n-1)^k<n^k$. Consequently every solution has at most $\lfloor k/\log2+\frac32\rfloor$ terms if $k<10^{55}$, and at most $k/\log2+\frac52$ terms in general.
\end{mainthm}

\begin{mainthm}[Theorem~\ref{thm:cubes}]\label{thmD}
For $k=3$, equation \eqref{eq:APintro} has solutions with exactly $n$ terms if and only if $n\in\{3,4,5\}$, and in each case it has infinitely many primitive solutions.
\end{mainthm}

\begin{mainthm}[Theorem~\ref{thm:even}, Proposition~\ref{prop:evencong}]\label{thmE}
Let $k$ be even. Then \eqref{eq:APintro} has no solution with three terms. If $k=2^s\ge4$, or if $4\le k\le24$, then \eqref{eq:APintro} has no solution at all.
\end{mainthm}

Theorem~\ref{thmC} replaces the exponential bound of \cite{DNZ} by a linear one, which cannot be improved by real considerations. The value $\lfloor k/\log2+\frac32\rfloor$ can be exceeded by one, but only for exponents attached to large partial quotients of $\log2$; the first of them has $330$ digits. Theorem~\ref{thmD} completes the treatment of cubes in \cite{DNZ}, where four to ten terms had remained open. The case $k=4$ of Theorem~\ref{thmE} is elementary and independent of Theorem~\ref{thmA}: it uses Lemma~\ref{lem:APreal}, the inequalities $1^4+\dots+6^4<7^4<8^4<1^4+\dots+7^4$, and congruences modulo $9$ and $16$. Congruences alone exclude only residue classes of the number of terms; the bound from the centre of mass makes the answer complete.

\subsection{The method}
The central idea of the paper is to read equations between power sums as balance conditions for suitable systems of masses. The use of mechanical ideas in mathematics is older than the calculus. Archimedes found the area of a parabolic segment and the volume of the sphere by balancing thin slices on a lever \cite{Heath}, and the barycentric calculus of M\"obius made the centre of mass a standard tool of geometry. In the twentieth century this tradition was kept alive in the books of Uspenskii \cite{Uspenskii}, of Balk and Boltyanskii \cite{BalkBoltyanskii}, and more recently of Levi \cite{Levi}, where many classical inequalities receive proofs by equilibrium arguments.

In \cite{Rahimi08,Rahimi09,Rahimi17} the second author proposed a way of bringing the centre of mass to bear on power sums. Given positive numbers $x_1<\dots<x_n$ and positive weights $a_1,\dots,a_n$, place at the point $x_i$ of the real line a mass $a_ix_i^{t-1}$. The centre of mass $c_t$ of this system is $c_t=S_t/S_{t-1}$, where $S_t=\sum a_ix_i^t$, and consequently
\begin{equation}\label{eq:product}
S_k=S_0\,c_1c_2\cdots c_k .
\end{equation}
This identity was used in \cite{Rahimi17} for equations of the form $ax^m+by^n=z^r$, and in \cite{MehdiNezhadRahimi} for Zeckendorf representations. The present paper adds a second mechanical ingredient, the moment of inertia, and uses two principles throughout.

The first is Steiner's parallel-axis theorem. The step from $c_t$ to $c_{t+1}$ is a Steiner step, $c_{t+1}=c_t+\rho_t^2/c_t$ (Theorem~\ref{thm:steiner}), and the same theorem gives the constant $\frac32$ and, in general, the first constant of each balance equation of Section~\ref{sec:type}. The second is the \emph{half-mass principle}. Equation \eqref{eq:EM} says that, in the system carrying the mass $i^k$ at each point $i=1,\dots,m$, the largest point carries one half of the total mass (Proposition~\ref{prop:EMreform}); the masses to its left decrease at least geometrically, and a geometric progression can balance its leading term only if its ratio exceeds $\frac12$, which is where $\log2$ comes from. For the progression equation the same principle shows that the most favourable configuration of $n$ terms is the Erd\H{o}s--Moser system on $0,1,\dots,n$; this is Theorem~\ref{thmC}, and Theorem~\ref{thmA} makes the resulting bound explicit.

A centre of mass is a ratio, and the balance condition determines an analytic parameter, here the exponent $\kappa_m$, which must be an integer for a solution to exist. The method is therefore suited to such \emph{balance equations}, for which it gives sharp necessary conditions. It cannot by itself decide whether a given integer is zero or a perfect power. Accordingly, the proofs of Theorems~\ref{thmA}--\ref{thmC} are analytic, and arithmetic enters only in Theorems~\ref{thmD} and~\ref{thmE}, through elliptic curves and congruences. The analytic proofs do not use the product formula \eqref{eq:product} itself; whether the asymptotic relation can be derived from it is Problem~\ref{prob:product}.

\subsection{Organisation}
Section~\ref{sec:systems} sets up the framework: the iterated mass systems, the Steiner step and Lagrange's identity. It also records, in mechanical form, some classical facts which follow at once: the monotonicity of Lehmer means, two-sided bounds for the moments of inertia, and the power-mean inequality. Section~\ref{sec:window} shows the half-mass principle in its simplest setting, equations $\sum a_ix_i^k=z^k$, where it explains the threshold of \cite{Rahimi09}. Section~\ref{sec:EM} reformulates \eqref{eq:EM} as a balance condition. Section~\ref{sec:threehalves} identifies the constant $\frac32$ and proves the second half of Theorem~\ref{thmA} with explicit constants, and Section~\ref{sec:uniform} proves the uniform estimate by means of Hermite's formula for the Hurwitz zeta function. Section~\ref{sec:type} treats general balance equations and proves Theorem~\ref{thmB}. Section~\ref{sec:AP} contains the applications to arithmetic progressions, Theorems~\ref{thmC}--\ref{thmE}. Some problems are collected in Section~\ref{sec:problems}, and Section~\ref{sec:earlier} corrects some statements of the earlier papers.

\section{Iterated mass systems}\label{sec:systems}

Lemma~\ref{lem:lever}, Theorem~\ref{thm:steiner} and Proposition~\ref{prop:lagrange} are used later; the other statements of this section are classical inequalities in mechanical form.

Throughout, $n\ge 2$, $0<x_1<x_2<\dots<x_n$ are real numbers and $a_1,\dots,a_n$ are positive weights. We write
\[
S_t=\sum_{i=1}^n a_ix_i^t\qquad(t\ge 0),
\]
so that $S_0=a_1+\dots+a_n$.

\begin{definition}\label{def:system}
For an integer $t\ge 1$, the \emph{$t$-th mass system} $\Sigma_t$ consists of the points $x_1,\dots,x_n$ of the real line, the point $x_i$ carrying the mass $m_i^{(t)}=a_ix_i^{t-1}$. Its total mass is $S_{t-1}$, and we denote by
\[
c_t=\frac{1}{S_{t-1}}\sum_{i=1}^n m_i^{(t)}x_i=\frac{S_t}{S_{t-1}}
\]
its centre of mass, by
\[
J_t=\sum_{i=1}^n m_i^{(t)}(x_i-c_t)^2
\]
its moment of inertia about the centre, by $\rho_t=(J_t/S_{t-1})^{1/2}$ its radius of gyration, and by
\[
p_t=\frac{a_nx_n^{t-1}}{S_{t-1}}
\]
the share of the total mass carried by the largest point $x_n$.
\end{definition}

The system $\Sigma_1$ carries the weights $a_i$ themselves; in the notation of \cite{Rahimi17}, $\Sigma_t$ is the $t$-mass-system with coefficients $a_1,\dots,a_n$. The centre $c_t=\sum a_ix_i^t/\sum a_ix_i^{t-1}$ is the weighted Lehmer mean of order $t$ of $x_1,\dots,x_n$ \cite{Lehmer}, and its monotonicity in $t$ is classical; the mechanical point of view adds the exact form of the increment in Theorem~\ref{thm:steiner}. The product formula \eqref{eq:product} follows by telescoping.

We record the precise form of the lever principle, which in the earlier papers \cite{Rahimi08,Rahimi17,MehdiNezhadRahimi} was used in an informal way.

\begin{lemma}[Lever principle]\label{lem:lever}
Let $\Sigma$ be a system of masses $\mu_1,\dots,\mu_n>0$ at points $x_1,\dots,x_n$, with centre of mass $c$, and let $q$ be a real number. Then $c\le q$ if and only if $\sum_i\mu_i(x_i-q)\le 0$, with equality on both sides simultaneously.
\end{lemma}

\begin{proof}
$\sum_i\mu_i(x_i-q)=\bigl(\sum_i\mu_i\bigr)(c-q)$.
\end{proof}

In words, the centre lies to the left of $q$ exactly when the total moment about $q$ is not positive. It is the balance of moments, and not the number of points on either side of $q$, that decides the position of the centre.

Some statements of the earlier papers are corrected in Section~\ref{sec:earlier}.

\subsection{The step from one centre to the next}

The moment of inertia of $\Sigma_t$ about the origin is
\[
\sum_{i=1}^n m_i^{(t)}x_i^2=\sum_{i=1}^n a_ix_i^{t+1}=S_{t+1}.
\]
Steiner's parallel-axis theorem, applied to $\Sigma_t$ and the origin, reads
\begin{equation}\label{eq:steiner}
S_{t+1}=J_t+S_{t-1}c_t^2 ;
\end{equation}
directly, expand $(x_i-c_t)^2$ and use $\sum_i m_i^{(t)}x_i=S_{t-1}c_t$.

\begin{theorem}\label{thm:steiner}
For every $t\ge 1$,
\[
c_{t+1}=c_t+\frac{\rho_t^2}{c_t}.
\]
\end{theorem}

\begin{proof}
Divide \eqref{eq:steiner} by $S_t=S_{t-1}c_t$. The left side becomes $c_{t+1}$, and the right side becomes $J_t/S_t+c_t=\rho_t^2/c_t+c_t$.
\end{proof}

Since $n\ge 2$, we have $\rho_t>0$; this is a sharper form of the monotonicity in \cite[Lemma~2.5]{Rahimi17}.

\begin{corollary}\label{cor:monotone}
The centres satisfy $x_1<c_1<c_2<c_3<\cdots<x_n$, and $c_t\to x_n$ as $t\to\infty$.
\end{corollary}

\begin{proof}
Strict increase follows from Theorem~\ref{thm:steiner}, and $x_1<c_t<x_n$ because every $\Sigma_t$ has positive mass at both $x_1$ and $x_n$. For the limit, Lemma~\ref{lem:lever} with $q=x_n$ gives
\[
x_n-c_t=\frac{1}{S_{t-1}}\sum_{i<n}a_ix_i^{t-1}(x_n-x_i)\le (x_n-x_1)(1-p_t),
\]
and $1-p_t\le a_n^{-1}S_0\,(x_{n-1}/x_n)^{t-1}$, which tends to $0$.
\end{proof}

\subsection{Two-sided estimates}

To estimate the increment in Theorem~\ref{thm:steiner} we bound the radius of gyration. The lower bound comes from Lagrange's identity, which expresses the moment of inertia about the centre through the mutual distances of the particles.

\begin{proposition}[Lagrange]\label{prop:lagrange}
For every $t\ge 1$,
\[
J_t=\frac{1}{S_{t-1}}\sum_{1\le i<j\le n}a_ia_j\,(x_ix_j)^{t-1}(x_j-x_i)^2 .
\]
Equivalently, $S_{t+1}S_{t-1}-S_t^2=\sum_{i<j}a_ia_j(x_ix_j)^{t-1}(x_j-x_i)^2$.
\end{proposition}

\begin{proof}
By \eqref{eq:steiner}, $S_{t-1}J_t=S_{t+1}S_{t-1}-S_t^2$. Now
\[
S_{t+1}S_{t-1}-S_t^2=\sum_{i,j}a_ia_j(x_ix_j)^{t-1}(x_i^2-x_ix_j)
=\frac12\sum_{i,j}a_ia_j(x_ix_j)^{t-1}(x_i-x_j)^2,
\]
the second equality being obtained by interchanging $i$ and $j$ and taking the mean of the two expressions.
\end{proof}

The upper bound is the inequality of Bhatia and Davis \cite{BhatiaDavis}, which has an equally simple mechanical proof.

\begin{lemma}\label{lem:BD}
Let $\Sigma$ be a system of positive masses on the interval $[\alpha,\beta]$, with centre $c$ and radius of gyration $\rho$. Then $\rho^2\le(\beta-c)(c-\alpha)$, with equality if and only if all the mass lies at the two end points.
\end{lemma}

\begin{proof}
Let $\mu_i$ be the masses, $x_i$ their positions and $\mu=\sum\mu_i$. Each term of $\sum_i\mu_i(\beta-x_i)(x_i-\alpha)$ is non-negative. Expanding, and using $\sum\mu_ix_i=\mu c$ and $\sum\mu_ix_i^2=\mu(\rho^2+c^2)$, which is \eqref{eq:steiner} in general form, we obtain
\[
0\le\sum_i\mu_i(\beta-x_i)(x_i-\alpha)=\mu\bigl((\alpha+\beta)c-\alpha\beta-\rho^2-c^2\bigr)=\mu\bigl((\beta-c)(c-\alpha)-\rho^2\bigr).
\]
Equality holds if and only if every $x_i$ is $\alpha$ or $\beta$.
\end{proof}

\begin{theorem}\label{thm:twosided}
Let $\delta=x_n-x_{n-1}$. For every $t\ge 1$,
\[
\delta^2\,p_t(1-p_t)\;\le\;\rho_t^2\;\le\;(x_n-c_t)(c_t-x_1),
\]
and consequently
\[
\frac{\delta^2\,p_t(1-p_t)}{c_t}\;\le\;c_{t+1}-c_t\;\le\;\frac{(x_n-c_t)(c_t-x_1)}{c_t}.
\]
\end{theorem}

\begin{proof}
The upper bound is Lemma~\ref{lem:BD} with $[\alpha,\beta]=[x_1,x_n]$. For the lower bound, keep in Proposition~\ref{prop:lagrange} only the pairs $(i,n)$ with $i<n$. Since $x_n-x_i\ge\delta$ for these pairs,
\[
S_{t-1}J_t\ge\delta^2\,a_nx_n^{t-1}\sum_{i<n}a_ix_i^{t-1}=\delta^2\,a_nx_n^{t-1}\bigl(S_{t-1}-a_nx_n^{t-1}\bigr).
\]
Dividing by $S_{t-1}^2$ gives $\rho_t^2\ge\delta^2p_t(1-p_t)$. The second chain follows from Theorem~\ref{thm:steiner}.
\end{proof}

\begin{remark}
The lower bound is largest when $p_t=\tfrac12$, that is, when the largest point carries exactly half of the mass of $\Sigma_t$. The centre then moves by at least $\delta^2/(4c_t)$. This is the first appearance of the half-mass condition which will be decisive in Sections~\ref{sec:window} and~\ref{sec:EM}.
\end{remark}

\subsection{Power means}

The product formula \eqref{eq:product} can be read as follows. Let
\[
M_k=\Bigl(\frac{S_k}{S_0}\Bigr)^{1/k}
\]
denote the weighted power mean of order $k$ of $x_1,\dots,x_n$. By \eqref{eq:product}, $M_k$ is the geometric mean of the first $k$ centres:
\[
M_k=(c_1c_2\cdots c_k)^{1/k}.
\]

\begin{corollary}\label{cor:powermean}
For every $k\ge 1$ we have $M_k\le c_k$, with equality only for $k=1$, and $M_1<M_2<M_3<\cdots$.
\end{corollary}

\begin{proof}
The geometric mean of the increasing numbers $c_1<\dots<c_k$ does not exceed the largest of them, and is strictly smaller when $k\ge2$. Hence
\[
M_{k+1}^{k+1}=M_k^k\,c_{k+1}>M_k^k\,c_k\ge M_k^{k+1}.\qedhere
\]
\end{proof}

Thus the power-mean inequality for integral orders is a consequence of Steiner's theorem.

\section{Windows for perfect powers}\label{sec:window}

This section illustrates the half-mass principle in its simplest setting. From now on the points $x_1<\dots<x_n$ and the weights $a_1,\dots,a_n$ are positive integers. We consider the equation
\begin{equation}\label{eq:general}
a_1x_1^k+a_2x_2^k+\dots+a_nx_n^k=z^k
\end{equation}
in a positive integer $z$, where $k\ge 2$. A solution $z$ must lie in a certain interval, which for large $k$ contains no integer.

\begin{lemma}\label{lem:log}
For $u>0$ we have $\log(1+u)>\dfrac{2u}{2+u}$.
\end{lemma}

\begin{proof}
The function $f(u)=\log(1+u)-2u/(2+u)$ vanishes at $u=0$ and $f'(u)=u^2/\bigl((1+u)(2+u)^2\bigr)>0$ for $u>0$.
\end{proof}

\begin{theorem}\label{thm:window}
If \eqref{eq:general} holds for a positive integer $z$, then
\[
x_n+1\le z<S_0^{1/k}\,c_k<S_0^{1/k}\,x_n ,
\]
where $c_k$ is the centre of mass of $\Sigma_k$. In particular
\[
k<\frac{\log S_0}{\log(1+1/x_n)}\qquad\text{and}\qquad x_n>\frac{k}{\log S_0}-\frac12 .
\]
\end{theorem}

\begin{proof}
Since $z^k\ge a_nx_n^k+a_1x_1^k>x_n^k$, we have $z\ge x_n+1$. By \eqref{eq:product} and Corollary~\ref{cor:monotone}, $z^k=S_0c_1\cdots c_k<S_0c_k^k$, and $c_k<x_n$. Hence $(1+1/x_n)^k\le(z/x_n)^k<S_0$, which is the first inequality. By Lemma~\ref{lem:log} with $u=1/x_n$ this gives $2k/(2x_n+1)<\log S_0$, which is the second.
\end{proof}

The first inequality is \cite[Theorem~3]{Rahimi09}, where it is obtained by the same argument with $\max_tc_t$ in place of $c_k$; the bound $z<S_0^{1/k}c_k$ is \cite[Theorem~2.6]{Rahimi17}. The second inequality is a slightly weaker but more transparent form of the first.

The case of two points with unit weights is the Fermat equation, for which $S_0=2$.

\begin{corollary}\label{cor:fermat}
Let $x<y$ and $k\ge 2$ be positive integers. If $x^k+y^k=z^k$ for a positive integer $z$, then
\[
y+1\le z<2^{1/k}\,c_k<2^{1/k}\,y,\qquad c_k=\frac{x^k+y^k}{x^{k-1}+y^{k-1}},
\]
and
\[
y>\frac{k}{\log 2}-\frac12 .
\]
\end{corollary}

By the theorem of Wiles \cite{Wiles}, Corollary~\ref{cor:fermat} is vacuous for $k\ge3$; it is included as the simplest illustration of Theorem~\ref{thm:window}.

\begin{remark}\label{rem:halfmass}
In mechanical terms the corollary says the following. In the system $\Sigma_{k+1}$ on the two points $x<y$, the larger point carries the share $y^k/z^k$ of the total mass $z^k$, and this share exceeds $\tfrac12$. The equation therefore requires $(z/y)^k<2$, while integrality requires $z/y\ge 1+1/y$; the two requirements are compatible only for $y>k/\log 2-\tfrac12$. In \cite[Remark~2]{Rahimi08} it was observed that for each fixed $y$ the equation has no solution for all sufficiently large $k$, and the explicit threshold $k\ge\log2/\log(1+1/y)$ was given in \cite[Corollary~4]{Rahimi09}; the half-mass interpretation explains why the constant $\log2$ appears.

For comparison, the classical argument $x^k=z^k-y^k\ge(y+1)^k-y^k>ky^{k-1}>kx^{k-1}$ shows that any solution satisfies $x>k$. The two bounds are of a different nature: the classical one concerns the smaller base, while the half-mass threshold concerns the larger, and it is the latter which generalizes to Theorem~\ref{thm:window} and to the equation of the next section.
\end{remark}

\section{The Erd\H{o}s--Moser equation}\label{sec:EM}

We reformulate \eqref{eq:EM} in terms of mass systems and show by a single mechanical argument that every solution satisfies $m>k/\log 2+\frac12$.

\begin{proposition}\label{prop:EMreform}
Let $m\ge2$ and $k\ge1$ be integers. The following are equivalent.
\begin{enumerate}[\rm(a)]
\item $1^k+2^k+\dots+(m-1)^k=m^k$.
\item In the system $\Sigma_{k+1}$ on the points $1,2,\dots,m$ with unit weights, the point $m$ carries exactly one half of the total mass.
\item $(m-1)\,c_1c_2\cdots c_k=m^k$, where $c_t$ is the centre of mass of the system $\Sigma_t$ on the points $1,2,\dots,m-1$ with unit weights.
\end{enumerate}
\end{proposition}

\begin{proof}
In (b) the point $i$ carries the mass $i^k$, so (b) says $m^k=\tfrac12\sum_{i=1}^m i^k$, which is (a). Condition (c) is (a) rewritten by means of \eqref{eq:product}, since here $S_0=m-1$.
\end{proof}

The masses $i^k$ in (b) decrease to the left of $m$ at least as fast as a geometric progression.

\begin{lemma}\label{lem:geometric}
Let $m\ge2$ and $k\ge1$, and put $\theta=(1-1/m)^k$, the ratio of the mass at $m-1$ to the mass at $m$ in the system of Proposition~\ref{prop:EMreform}(b). Then for $1\le j\le m-1$,
\[
(m-j)^k\le\theta^j\,m^k ,
\]
with strict inequality when $j\ge 2$.
\end{lemma}

\begin{proof}
By Bernoulli's inequality, $(1-1/m)^j\ge 1-j/m$, with strict inequality for $j\ge2$. Raising to the power $k$ gives $(1-j/m)^k\le\theta^j$.
\end{proof}

\begin{theorem}\label{thm:EM}
If $m\ge2$ and $k\ge1$ satisfy $1^k+2^k+\dots+(m-1)^k=m^k$, then $m\ge3$ and
\[
m>\frac{k}{\log 2}+\frac12 .
\]
\end{theorem}

\begin{proof}
For $m=2$ the equation reads $1=2^k$, which is impossible; hence $m\ge3$. By Lemma~\ref{lem:geometric}, the total mass to the left of $m$ in the system of Proposition~\ref{prop:EMreform}(b) satisfies
\[
m^k=\sum_{j=1}^{m-1}(m-j)^k<m^k\sum_{j=1}^{m-1}\theta^j<m^k\,\frac{\theta}{1-\theta},
\]
the first inequality being strict because the term $j=2$ is present. Hence $\theta>\frac12$, that is,
\[
k\log\Bigl(1+\frac{1}{m-1}\Bigr)=-k\log\Bigl(1-\frac1m\Bigr)<\log 2 .
\]
By Lemma~\ref{lem:log} with $u=1/(m-1)$, the left side exceeds $2k/(2m-1)$. Therefore $2m-1>2k/\log 2$, which is the assertion.
\end{proof}

\begin{remark}
By \cite{GMZ}, large solutions satisfy $m=k/\log2+\frac32+O(1/m)$. Thus Theorem~\ref{thm:EM} gives the main term exactly and the constant term up to one unit; the correct constant is obtained in Section~\ref{sec:threehalves}.
\end{remark}

\begin{remark}\label{rem:heuristic}
The product formula in Proposition~\ref{prop:EMreform}(c) shows where the constant $\log 2$ comes from. Write it as
\[
\sum_{t=1}^k\log\frac{m}{c_t}=\log(m-1).
\]
For $1\ll t\le k$ the masses $i^{t-1}$ near the point $m-1$ behave like a geometric progression with ratio $e^{-(t-1)/m}$ in the distance from $m$, and the centre of such a progression lies at distance about $\bigl(1-e^{-t/m}\bigr)^{-1}$ to the left of $m$. Replacing the sum by an integral with $s=t/m$, the left side becomes approximately
\[
\int_{1/m}^{k/m}\frac{ds}{1-e^{-s}}=\log\bigl(e^{k/m}-1\bigr)-\log\bigl(e^{1/m}-1\bigr)\approx\log\bigl(e^{k/m}-1\bigr)+\log m ,
\]
and the equation becomes $e^{k/m}-1\approx1$, that is, $k/m\approx\log2$. We have not made this argument rigorous (Problem~\ref{prob:product}); the next section follows a different route.
\end{remark}

\section{The constant $\tfrac32$}\label{sec:threehalves}

Throughout this section we write $L=\log 2$ and
\[
\gamma_0=\frac{25}{12}-3\log 2=0.003891\ldots .
\]

\subsection{The half-geometric system}

\begin{definition}\label{def:gamma}
Fix a point $P$ of the real line. For $\lambda>0$ let $\Gamma_\lambda$ be the system which carries the mass $e^{-\lambda j}$ at the point $P-j$, for $j=1,2,3,\dots$. Its moments about $P$ are
\[
Q_r(\lambda)=\sum_{j=1}^\infty j^r e^{-\lambda j}\qquad(r\ge 0).
\]
Thus $Q_0$ is the total mass of $\Gamma_\lambda$, $Q_1$ its static moment and $Q_2$ its moment of inertia about $P$. We write $\Gamma=\Gamma_L$; it carries the masses $\tfrac12,\tfrac14,\tfrac18,\dots$.
\end{definition}

The series converge uniformly on every half-line $\lambda\ge\lambda_0>0$, together with their derivatives, so that
\begin{equation}\label{eq:Qderiv}
Q_r'(\lambda)=-Q_{r+1}(\lambda),
\end{equation}
and each $Q_r$ is a decreasing function of $\lambda$.

\begin{lemma}\label{lem:Qvalues}
We have $Q_0(L)=1$ and, for $r\ge1$,
\[
Q_r(L)=1+\sum_{i=0}^{r-1}\binom{r}{i}Q_i(L).
\]
In particular $Q_0(L),Q_1(L),\dots,Q_7(L)$ are $1,2,6,26,150,1082,9366,94586$.
\end{lemma}

\begin{proof}
Write $Q_r=Q_r(L)=\sum_{j\ge1}j^r2^{-j}$. Replacing $j$ by $j+1$,
\[
Q_r=\frac12\sum_{j\ge0}(j+1)^r2^{-j}=\frac12\Bigl(1+\sum_{i=0}^{r}\binom ri Q_i\Bigr),
\]
where the term $1$ comes from $j=0$ and the index $i=0$. Solving for $Q_r$ gives the recurrence, and $Q_0=1$ is the sum of the geometric series. The numerical values follow successively.
\end{proof}

For $r\ge1$ the numbers $\frac12Q_r(L)$ are the ordered Bell numbers, but we shall need only the values listed.

\begin{proposition}\label{prop:Gamma}
The system $\Gamma$ has total mass $1$, its centre of mass lies at the distance $d=2$ from $P$, and its radius of gyration about the centre is $\rho=\sqrt2$. Consequently
\[
\frac{Q_2(L)}{2Q_1(L)}=\frac{\rho^2+d^2}{2d}=\frac32 .
\]
\end{proposition}

\begin{proof}
By Lemma~\ref{lem:Qvalues}, the total mass is $Q_0=1$ and $d=Q_1/Q_0=2$. By Steiner's theorem the moment of inertia about $P$ is $Q_2=Q_0(\rho^2+d^2)$, so $\rho^2=6-4=2$. Finally $Q_1=Q_0d$, whence $Q_2/2Q_1=(\rho^2+d^2)/2d=6/4$.
\end{proof}

Since $Q_0(L)=1$, the system $\Gamma$ balances, in mass, a unit mass placed at $P$. In the language of Proposition~\ref{prop:EMreform}(b), $\Gamma$ together with a unit mass at $P$ is the limiting form of a system whose largest point carries exactly one half of the total mass.

We shall need the moments $Q_r(\lambda)$ for $\lambda$ slightly smaller than $L$.

\begin{lemma}\label{lem:inflation}
If $0\le a\le 1.1\cdot10^{-5}$ and $1\le r\le 7$, then $Q_r(L)\le Q_r(L-a)\le 1.0002\,Q_r(L)$.
\end{lemma}

\begin{proof}
The first inequality holds because $Q_r$ is decreasing. For the second we use Euler's formula
\[
\sum_{j=1}^\infty j^rq^j=\frac{q\,A_r(q)}{(1-q)^{r+1}}\qquad(0<q<1),
\]
where $A_r$ is the Eulerian polynomial of degree $r-1$, whose coefficients are non-negative integers. Put $t=e^a$, so that $e^{-(L-a)}=t/2$. Since $qA_r(q)$ is a polynomial with non-negative coefficients in which only the powers $q,q^2,\dots,q^r$ occur, replacing $q=\frac12$ by $q=\frac t2$ multiplies it by at most $t^r$; and $1-\frac t2=\frac12(2-t)$. Hence
\[
Q_r(L-a)\le\frac{t^r}{(2-t)^{r+1}}\,Q_r(L)\le\frac{t^7}{(2-t)^8}\,Q_r(L).
\]
Now $t\le1+1.11\cdot10^{-5}$, so $t^7\le1.000078$ and $(2-t)^8\ge0.99991$, and the quotient is less than $1.0002$.
\end{proof}

\subsection{The Erd\H{o}s--Moser system seen from its largest point}

For an integer $m\ge3$ and a real number $\kappa>0$ put
\[
F_m(\kappa)=\sum_{j=1}^{m-1}\Bigl(1-\frac jm\Bigr)^{\kappa}.
\]
This is the total mass to the left of the point $m$ in the system which carries the mass $i^\kappa$ at each point $i=1,\dots,m$, measured in units of the mass at $m$. The function $F_m$ decreases strictly from $F_m(0)=m-1\ge2$ to $0$, so there is a unique $\kappa_m>0$ with $F_m(\kappa_m)=1$. By Proposition~\ref{prop:EMreform}, a pair $(m,k)$ satisfies \eqref{eq:EM} if and only if $\kappa_m=k$.

Let $\lambda=\kappa/m$. The mass at distance $j$ from the point $m$ is
\[
\Bigl(1-\frac jm\Bigr)^{\kappa}=e^{-\lambda j}\exp\Bigl(\kappa\Bigl(\log\Bigl(1-\frac jm\Bigr)+\frac jm\Bigr)\Bigr),
\]
and the exponent in the second factor is approximately $-\lambda j^2/2m$. Thus the system is obtained from $\Gamma_\lambda$, placed with $P=m$, by diminishing the $j$-th mass by a factor close to $1-\lambda j^2/2m$, that is, in proportion to $j^2$, which is the share of that mass in the moment of inertia about $P$. To a first approximation the total mass to the left of $P$ is therefore $Q_0(\lambda)-(\lambda/2m)Q_2(\lambda)$, and the balance requires this to equal $1$. At $\lambda=L$ there is a deficit of about $(L/2m)Q_2(L)$; lowering $\lambda$ from $L$ to $L-a$ raises $Q_0$ by about $Q_1(L)\,a$, by \eqref{eq:Qderiv}. Hence $a\approx LQ_2(L)/(2mQ_1(L))$, and by Proposition~\ref{prop:Gamma}
\[
\kappa_m=m(L-a)\approx mL-\tfrac32L,\qquad\text{that is,}\qquad m\approx\frac{\kappa_m}{\log 2}+\frac32 .
\]
The following theorem makes this precise and carries it one step further.

\begin{theorem}\label{thm:main}
Let $m\ge10^5$ and define $\gamma_m$ by
\[
\kappa_m=\log 2\,\Bigl(m-\frac32-\frac{\gamma_m}{m}\Bigr),\qquad\text{equivalently}\qquad \frac{\kappa_m}{m}=\log2\,\Bigl(1-\frac{3}{2m}-\frac{\gamma_m}{m^2}\Bigr).
\]
Then
\[
-\frac{87}{m}<\gamma_m-\gamma_0<\frac{116}{m}.
\]
\end{theorem}

\begin{corollary}\label{cor:convergent}
Let $m\ge10^5$ and $k$ be positive integers satisfying \eqref{eq:EM}. Then
\[
m=\frac{k}{\log 2}+\frac32+\frac{\gamma_m}{m},\qquad 0.003<\gamma_m<0.0051,
\]
and
\[
0<\log 2-\frac{2k}{2m-3}<\frac{0.0142}{(2m-3)^2}.
\]
Consequently $2k/(2m-3)$ is a convergent of even index of the continued fraction of $\log 2$.
\end{corollary}

\begin{proof}
Here $\kappa_m=k$, and the bounds for $\gamma_m$ follow from Theorem~\ref{thm:main}, since $\gamma_0-87\cdot10^{-5}>0.003$ and $\gamma_0+116\cdot10^{-5}<0.0051$. From the definition of $\gamma_m$,
\[
\log 2-\frac{2k}{2m-3}=\frac{2\gamma_m\log 2}{m(2m-3)}=\frac{2\gamma_m\log2}{(2m-3)^2}\cdot\frac{2m-3}{m}<\frac{4\gamma_m\log 2}{(2m-3)^2},
\]
and $4\cdot0.0051\cdot\log2<0.0142$. The right side is less than $1/(2(2m-3)^2)$; this inequality persists when the fraction $2k/(2m-3)$ is reduced to lowest terms, since its denominator can only decrease. By the theorem of Legendre \cite[Theorem~184]{HW}, $2k/(2m-3)$ is a convergent of $\log2$, and since it is smaller than $\log2$ its index is even \cite[Chapter~X]{HW}.
\end{proof}

\begin{remark}\label{rem:GMZ}
The asymptotic relation $\kappa_m=\log2\,(m-\frac32-\gamma_0/m+O(m^{-2}))$ was obtained in \cite[Theorem~1]{GMZ}, together with $0<\gamma_m<0.004$ for $m>10^9$, and Corollary~\ref{cor:convergent} recovers \cite[Corollary~1]{GMZ}; combined with Moser's arithmetic conditions \cite{Moser}, the latter is the starting point of the bound $m>2.7\cdot10^{10^9}$ for solutions other than $1+2=3$. The proof given here is different and self-contained: it gives explicit constants from $m=10^5$ on, uses no arithmetic of the equation, and identifies the constant $\frac32$ mechanically (Proposition~\ref{prop:Gamma}).
\end{remark}

\subsection{Proof of Theorem~\ref{thm:main}}

For $1\le j\le m-1$ and $\lambda=\kappa/m$ we write
\[
w_j=\Bigl(1-\frac jm\Bigr)^{\kappa},\qquad x_j=\frac{\lambda j^2}{2m},\qquad y_j=\frac{\lambda j^3}{3m^2},\qquad v_j=\frac{\lambda j^4}{2m^3}.
\]

\begin{lemma}\label{lem:termwise}
\begin{enumerate}[\rm(a)]
\item For $1\le j\le m-1$,
\[
w_j\le e^{-\lambda j}\Bigl(1-x_j-y_j+\tfrac12(x_j+y_j)^2\Bigr).
\]
\item For $1\le j\le m/2$,
\[
w_j\ge e^{-\lambda j}\Bigl(1-x_j-y_j-v_j+\tfrac12x_j^2-\tfrac16x_j^3\Bigl(1+\frac{6j}{m}\Bigr)\Bigr).
\]
\end{enumerate}
\end{lemma}

\begin{proof}
Let $u=j/m$. From $-\log(1-u)=\sum_{r\ge1}u^r/r$ we have
\[
-\log(1-u)\ge u+\frac{u^2}2+\frac{u^3}3,
\]
and, when $u\le\frac12$,
\[
-\log(1-u)\le u+\frac{u^2}2+\frac{u^3}3+\frac{u^4}{4(1-u)}\le u+\frac{u^2}2+\frac{u^3}3+\frac{u^4}2 .
\]
Multiplying by $\kappa=\lambda m$, we obtain $w_j\le e^{-\lambda j}e^{-(x_j+y_j)}$ for all $j$, and $w_j\ge e^{-\lambda j}e^{-z_j}$ with $z_j=x_j+y_j+v_j$ for $j\le m/2$. We now use the inequalities
\[
1-s+\frac{s^2}2-\frac{s^3}6\le e^{-s}\le1-s+\frac{s^2}2\qquad(s\ge0).
\]
The upper one gives (a). For (b), note that $z_j^2\ge x_j^2$ and $z_j=x_j(1+\sigma)$ with $\sigma=2j/3m+j^2/m^2$. For $j\le m/2$ we have $\sigma\le\frac76\cdot\frac jm\le\frac7{12}$, so that
\[
(1+\sigma)^3=1+\sigma(3+3\sigma+\sigma^2)\le1+5.1\,\sigma\le1+\frac{6j}{m}.
\]
Hence $z_j^3\le x_j^3(1+6j/m)$, and (b) follows from the lower inequality.
\end{proof}

Put
\[
\psi(\lambda)=\frac{\lambda^2}8Q_4(\lambda)-\frac\lambda3Q_3(\lambda),\qquad
\Phi(\lambda)=Q_0(\lambda)-\frac{\lambda}{2m}Q_2(\lambda)+\frac{\psi(\lambda)}{m^2}.
\]

\begin{lemma}\label{lem:FPhi}
Let $m\ge10^5$ and $L-1.1\cdot10^{-5}\le\lambda\le L$, and put $\kappa=\lambda m$. Then
\[
\Phi(\lambda)-\frac{118}{m^3}\le F_m(\kappa)\le\Phi(\lambda)+\frac{87}{m^3}.
\]
\end{lemma}

\begin{proof}
The bracket in Lemma~\ref{lem:termwise}(a) is positive, because $1-s+s^2/2>0$ for every real $s$. Summing (a) over $1\le j\le m-1$ and then extending the sum to all $j\ge1$ we obtain, since $\frac12(x_j+y_j)^2=\lambda^2j^4/8m^2+\lambda^2j^5/6m^3+\lambda^2j^6/18m^4$,
\[
F_m(\kappa)\le Q_0-\frac{\lambda}{2m}Q_2-\frac{\lambda}{3m^2}Q_3+\frac{\lambda^2}{8m^2}Q_4+\frac{\lambda^2}{6m^3}Q_5+\frac{\lambda^2}{18m^4}Q_6
=\Phi+\frac{\lambda^2}{6m^3}Q_5+\frac{\lambda^2}{18m^4}Q_6,
\]
all moments being taken at $\lambda$. By Lemmas~\ref{lem:Qvalues} and~\ref{lem:inflation} and $\lambda\le L$, we have $\lambda^2Q_5/6\le L^2\cdot1082.3/6<86.67$ and $\lambda^2Q_6/18m<0.01$, which gives the upper bound.

For the lower bound let $N=\lfloor m/2\rfloor$. We sum (b) over $1\le j\le N$ and discard the remaining terms $w_j$, which are positive. In the sums which carry a negative sign we may extend the summation to all $j\ge1$, since this only decreases the right side. Writing $T_r=\sum_{j>N}j^re^{-\lambda j}$, we obtain
\[
F_m(\kappa)\ge\Phi-T_0-\frac{\lambda^2}{8m^2}T_4-\frac{\lambda}{2m^3}Q_4-\frac{\lambda^3}{48m^3}Q_6-\frac{\lambda^3}{8m^4}Q_7 .
\]
By Lemma~\ref{lem:inflation}, $\lambda Q_4/2\le L\cdot150.03/2<52.0$, $\lambda^3Q_6/48\le L^3\cdot9369.9/48<65.1$ and $\lambda^3Q_7/8m<0.04$. As for the tails, $N\ge(m-1)/2>200$, and for $j\ge200$ we have $j^4\le2^{j/6}$; moreover $e^{-\lambda j}\le2^{-j}e^{1.1\cdot10^{-5}j}\le2^{-5j/6}$. Hence
\[
T_0+T_4\le2\sum_{j>N}2^{-2j/3}<6\cdot2^{-m/3},
\]
which is far smaller than $m^{-4}$. Altogether the subtracted terms are less than $118/m^3$.
\end{proof}

\begin{lemma}\label{lem:Phi}
Let $m\ge10^5$, let $|c|\le1$, and put $\lambda=L(1-3/2m-c/m^2)$. Then
\[
-\frac{42}{m^3}\le\Phi(\lambda)-1-\frac{2L(c-\gamma_0)}{m^2}\le\frac{33}{m^3}.
\]
\end{lemma}

\begin{proof}
Put $a=L-\lambda=A/m+B/m^2$, where $A=\frac32L$ and $B=cL$. Then $0<a\le(A+L/m)/m<1.04/m\le1.04\cdot10^{-5}$, so that Lemma~\ref{lem:inflation} applies on the interval $[L-a,L]$. By Taylor's formula with Lagrange's remainder, \eqref{eq:Qderiv} and Lemma~\ref{lem:Qvalues},
\begin{align*}
Q_0(L-a)&=1+2a+3a^2+\varrho_0, & 0&\le\varrho_0=\tfrac{a^3}6Q_3(L-\xi_0)\le4.34\,a^3,\\
Q_2(L-a)&=6+26a+\varrho_2, & 0&\le\varrho_2=\tfrac{a^2}2Q_4(L-\xi_2)\le75.02\,a^2,
\end{align*}
with $0<\xi_0,\xi_2<a$. Next, $\psi(L)=\frac{75}4L^2-\frac{26}3L$, and
\[
\psi'(\lambda)=\frac{7\lambda}{12}Q_4(\lambda)-\frac{\lambda^2}8Q_5(\lambda)-\frac13Q_3(\lambda).
\]
For $L-a\le\lambda\le L$, Lemma~\ref{lem:inflation} shows that the first term lies between $60.64$ and $60.67$, and that the sum of the last two lies between $73.64$ and $73.67$. Hence $|\psi'|<13.1$ on $[L-a,L]$, and $\psi(L-a)=\psi(L)+\varrho_3$ with $|\varrho_3|\le13.1\,a$.

Since $-\frac{\lambda}{2m}Q_2(\lambda)=-\frac{L-a}{2m}(6+26a+\varrho_2)$, we find
\[
\Phi(\lambda)-1=2a+3a^2-\frac{3L}{m}-\frac{(13L-3)a}{m}+\frac{13a^2}{m}+\frac{\psi(L)}{m^2}+\varrho_0-\frac{(L-a)\varrho_2}{2m}+\frac{\varrho_3}{m^2}.
\]
We insert $a=A/m+B/m^2$. The terms of order $1/m$ are $2A/m-3L/m=0$; this cancellation is the relation $A=LQ_2(L)/2Q_1(L)$ of the heuristic argument above. The terms of order $1/m^2$ are
\[
\frac1{m^2}\Bigl(2B+3A^2-(13L-3)A+\psi(L)\Bigr)=\frac1{m^2}\Bigl(2cL+6L^2-\frac{25}6L\Bigr)=\frac{2L(c-\gamma_0)}{m^2}.
\]
The remaining terms are
\[
\frac{B(6A-13L+3)}{m^3}+\frac{3B^2}{m^4}+\frac{13a^2}{m}+\varrho_0-\frac{(L-a)\varrho_2}{2m}+\frac{\varrho_3}{m^2}.
\]
Here $B(6A-13L+3)=cL(3-4L)$, whose absolute value is less than $0.158$; further $3B^2/m^4<10^{-4}/m^3$, $0\le13a^2/m<14.07/m^3$, $0\le\varrho_0<4.89/m^3$, $0\le(L-a)\varrho_2/2m<28.13/m^3$ and $|\varrho_3|/m^2<13.63/m^3$. Adding the positive contributions gives the upper bound $32.75/m^3$, and adding the negative ones gives the lower bound $-41.92/m^3$.
\end{proof}

\begin{proof}[Proof of Theorem~\ref{thm:main}]
For $|c|\le1$ put $\kappa(c)=mL(1-3/2m-c/m^2)=L(m-\frac32-c/m)$. Then $\lambda=\kappa(c)/m$ satisfies the hypothesis of Lemma~\ref{lem:FPhi}, by the bound for $a$ in the proof of Lemma~\ref{lem:Phi}. Combining the two lemmas,
\begin{equation}\label{eq:star}
-\frac{160}{m^3}\le F_m(\kappa(c))-1-\frac{2L(c-\gamma_0)}{m^2}\le\frac{120}{m^3}\qquad(|c|\le1).
\end{equation}
For $c=1$ the left inequality gives $F_m(\kappa(1))-1\ge(2L(1-\gamma_0)-160/m)/m^2>0$, and for $c=-1$ the right inequality gives $F_m(\kappa(-1))-1\le(-2L(1+\gamma_0)+120/m)/m^2<0$. Since $c\mapsto F_m(\kappa(c))$ is continuous and strictly increasing, the unique $c$ with $F_m(\kappa(c))=1$, which is $\gamma_m$, satisfies $|\gamma_m|<1$. Taking $c=\gamma_m$ in \eqref{eq:star} we obtain
\[
-\frac{120}{2Lm}\le\gamma_m-\gamma_0\le\frac{160}{2Lm},
\]
and $120/2L<86.6$, $160/2L<115.5$.
\end{proof}

\section{A uniform estimate}\label{sec:uniform}

For small $m$ the method of Section~\ref{sec:threehalves} loses precision, since it expands the factor $e^{-\lambda j^2/2m}$ in powers of $j^2/m$, and the resulting moments $Q_r(L)$ grow like $r!/L^{r+1}$. We use instead Hermite's formula for the Hurwitz zeta function, which writes $F_m(\kappa)$ as a function of $h=1/m$ that is analytic in a large disc, plus small remainders. As before $L=\log2$, and we write
\[
\beta=\frac{1}{8\log 2}=0.18033\ldots,\qquad h=\frac1m .
\]

\subsection{Hermite's formula}

Let
\[
u(z)=\frac{1}{e^z-1}-\frac1z+\frac12=\sum_{k=1}^\infty\frac{B_{2k}}{(2k)!}\,z^{2k-1}\qquad(|z|<2\pi),
\]
where $B_{2k}$ are the Bernoulli numbers. For real $\lambda>0$ we have $u(\lambda)=Q_0(\lambda)-1/\lambda+\frac12$, so that $u$ differs from the total mass of the system $\Gamma_\lambda$ of Definition~\ref{def:gamma} only by an elementary term; in particular $u(L)=\frac32-\frac1L$ and, for $k\ge1$,
\begin{equation}\label{eq:uk}
u_k:=u^{(k)}(L)=(-1)^k\Bigl(Q_k(L)-\frac{k!}{L^{k+1}}\Bigr).
\end{equation}
The numbers $u_k$ are small; for instance $u_1=0.0813\ldots$, $u_2=-0.0055\ldots$, $u_3=-0.0074\ldots$. This reflects the fact that $u$ is analytic in the disc $|z|<2\pi$, whereas $Q_0(\lambda)$ has a pole at $\lambda=0$. We shall also use the majorant
\[
\hat u(x)=\sum_{k=1}^\infty\frac{|B_{2k}|}{(2k)!}\,x^{2k-1}=\frac1x\Bigl(1-\frac x2\cot\frac x2\Bigr)\qquad(0<x<2\pi),
\]
which has non-negative coefficients and satisfies $|u^{(k)}(z)|\le\hat u^{(k)}(|z|)$ for $|z|<2\pi$ and $k\ge0$.

Two classical formulas will be used. The first is
\begin{equation}\label{eq:sinint}
2\int_0^\infty\frac{\sin\lambda t}{e^{2\pi t}-1}\,dt=u(\lambda)\qquad(\lambda>0),
\end{equation}
see \cite[3.911.2]{GR}. Differentiating with respect to $\lambda$ under the integral sign we obtain, for $j\ge0$,
\begin{equation}\label{eq:momints}
2\int_0^\infty\frac{t^{2j}\sin\lambda t}{e^{2\pi t}-1}\,dt=(-1)^ju^{(2j)}(\lambda),\qquad
2\int_0^\infty\frac{t^{2j+1}\cos\lambda t}{e^{2\pi t}-1}\,dt=(-1)^ju^{(2j+1)}(\lambda).
\end{equation}
The second is Hermite's formula for the Hurwitz zeta function $\zeta(s,a)=\sum_{n\ge0}(n+a)^{-s}$, valid for all complex $s\ne1$ and $a>0$ \cite[\S13.2]{WW}:
\[
\zeta(s,a)=\frac{a^{-s}}2+\frac{a^{1-s}}{s-1}+2\int_0^\infty\frac{\sin(s\arctan(t/a))}{(a^2+t^2)^{s/2}}\,\frac{dt}{e^{2\pi t}-1}.
\]

\begin{lemma}\label{lem:hermite}
For an integer $m\ge2$ and a real number $\kappa>0$,
\[
F_m(\kappa)=\frac{m}{\kappa+1}-\frac12+\frac{\zeta(-\kappa)}{m^\kappa}+2\int_0^\infty\operatorname{Im}\Bigl(1+\frac{it}m\Bigr)^{\kappa}\frac{dt}{e^{2\pi t}-1}.
\]
\end{lemma}

\begin{proof}
For $\operatorname{Re}s>1$ we have $\zeta(s)-\zeta(s,m)=\sum_{i=1}^{m-1}i^{-s}$, and by analytic continuation this holds for all $s\ne1$. Take $s=-\kappa$ and $a=m$ in Hermite's formula, observe that $(m^2+t^2)^{\kappa/2}\sin(\kappa\arctan(t/m))=m^\kappa\operatorname{Im}(1+it/m)^\kappa$, and divide by $m^\kappa$.
\end{proof}

Now put $\kappa=\lambda m$, so that $m/(\kappa+1)=1/(\lambda+h)$. By \eqref{eq:sinint}, Lemma~\ref{lem:hermite} becomes
\begin{equation}\label{eq:exact}
F_m(\kappa)-1=u(\lambda)-\frac32+\frac1{\lambda+h}+J(\lambda,h)+\frac{\zeta(-\kappa)}{m^\kappa},
\end{equation}
where
\[
J(\lambda,h)=2\int_0^\infty\operatorname{Im}\Bigl((1+iht)^{\lambda/h}-e^{i\lambda t}\Bigr)\frac{dt}{e^{2\pi t}-1}.
\]
The three terms $u(\lambda)-\frac32+1/(\lambda+h)$ vanish together at $h=0$, $\lambda=L$, which is once more the half-mass condition $Q_0(L)=1$.

\subsection{The integral $J$}

We have $(1+iht)^{\lambda/h}=e^{i\lambda t}e^{W}$ with $W=(\lambda/h)\bigl(\log(1+iht)-iht\bigr)$, and
\[
e^W=1+h\,\frac{\lambda t^2}2+h^2\Bigl(\frac{\lambda^2t^4}8-\frac{i\lambda t^3}3\Bigr)+h^3\Bigl(\frac{\lambda^3t^6}{48}-\frac{\lambda t^4}4-\frac{i\lambda^2t^5}6\Bigr)+\cdots .
\]
Inserting the terms up to $h^3$ into $J$ and using \eqref{eq:momints} we obtain $hP_1(\lambda)+h^2P_2(\lambda)+h^3P_3(\lambda)$, where
\begin{align*}
P_1(\lambda)&=-\frac{\lambda}{2}u''(\lambda),\qquad P_2(\lambda)=\frac{\lambda}{3}u'''(\lambda)+\frac{\lambda^2}{8}u^{(4)}(\lambda),\\
P_3(\lambda)&=-\frac{\lambda}{4}u^{(4)}(\lambda)-\frac{\lambda^2}{6}u^{(5)}(\lambda)-\frac{\lambda^3}{48}u^{(6)}(\lambda).
\end{align*}

\begin{lemma}\label{lem:Jrem}
Let $m\ge10$ and $0<\lambda\le L$. Then
\[
\bigl|J(\lambda,h)-hP_1(\lambda)-h^2P_2(\lambda)-h^3P_3(\lambda)\bigr|\le0.08\,h^4 .
\]
\end{lemma}

\begin{proof}
For fixed $t>0$ the function $\eta\mapsto W(\eta)=(\lambda/\eta)(\log(1+i\eta t)-i\eta t)$ is analytic in $|\eta|<1/t$. On the circle $|\eta|=3/4t$ we have $|\log(1+z)-z|\le\sum_{k\ge2}|z|^k/k\le\frac98$ with $|z|=\frac34$, hence $|W|\le\frac43\lambda t\cdot\frac98=\frac32\lambda t$. By Cauchy's estimates the coefficient $w_n(t)$ of $\eta^n$ in $e^{W}$ satisfies $|w_n(t)|\le e^{3\lambda t/2}(4t/3)^n$, and for $4ht/3\le\frac12$
\[
\Bigl|e^{W(h)}-\sum_{n=0}^3w_n(t)h^n\Bigr|\le e^{3\lambda t/2}\sum_{n\ge4}\Bigl(\frac{4ht}3\Bigr)^n\le2e^{3\lambda t/2}\Bigl(\frac{4ht}3\Bigr)^4 .
\]
The quantity to be estimated equals $2\int_0^\infty\operatorname{Im}\bigl(e^{i\lambda t}(e^{W}-\sum_{n\le3}w_nh^n)\bigr)\,dt/(e^{2\pi t}-1)$. We split the integral at $T=3/(8h)=3m/8\ge3.75$.

For $t\le T$ the bound above gives the contribution
\[
\le4\Bigl(\frac{4h}3\Bigr)^4\int_0^\infty\frac{t^4e^{bt}}{e^{2\pi t}-1}\,dt=4\Bigl(\frac{4h}3\Bigr)^4\cdot24\sum_{j\ge1}\frac1{(2\pi j-b)^5},\qquad b=\frac32\lambda\le\frac32L .
\]
Since $2\pi j-b\ge2\pi j(1-b/2\pi)$, the sum is at most $(1-3L/4\pi)^{-5}\zeta(5)(2\pi)^{-5}<2.471\cdot1.0370/9792<2.617\cdot10^{-4}$, and the contribution is less than $0.0795\,h^4$.

For $t>T$ we use $|(1+iht)^{\kappa}|=(1+h^2t^2)^{\kappa/2}\le(1+ht)^\kappa\le e^{\lambda t}$ and $|\sum_{n\le3}w_nh^n|\le e^{3\lambda t/2}(1+4ht/3)^3\le e^{3\lambda t/2}e^{4ht}$. With $h\le\frac1{10}$ both are at most $e^{(3L/2+0.4)t}$, and the contribution is at most
\[
4.001\int_T^\infty e^{-(2\pi-3L/2-0.4)t}\,dt<0.83\,e^{-1.81m}<0.0002\,h^4
\]
for $m\ge10$, since $m^4e^{-1.81m}$ decreases for $m\ge3$. Adding the two contributions gives the lemma.
\end{proof}

\subsection{The analytic part}

For a real parameter $c$ put $\lambda(h)=L(1-\frac32h-ch^2)$ and
\[
\Psi_c(h)=u(\lambda)-u(L)+\frac1{\lambda+h}-\frac1L+hP_1(\lambda)+h^2P_2(\lambda)+h^3P_3(\lambda),\qquad\lambda=\lambda(h).
\]
Since $u(L)-\frac32+\frac1L=0$, \eqref{eq:exact} and Lemma~\ref{lem:Jrem} give, for real $h=1/m$ with $m\ge10$ and $\kappa=\lambda(h)m$,
\begin{equation}\label{eq:FPsi}
\Bigl|F_m(\kappa)-1-\Psi_c(h)-\frac{\zeta(-\kappa)}{m^\kappa}\Bigr|\le0.08\,h^4 .
\end{equation}
The function $\Psi_c$ is defined by the same formula for complex $h$ and is analytic as long as $|\lambda(h)|<2\pi$ and $\lambda(h)+h\ne0$.

\begin{lemma}\label{lem:Psicoef}
The Taylor expansion of $\Psi_c$ at $h=0$ has the form $\Psi_c(h)=A_2(c)h^2+A_3(c)h^3+O(h^4)$, where
\[
A_2(c)=\frac L6\,(36L+12c-25)=2L(c-\gamma_0),\qquad
A_3(c)=-\frac L6\,\bigl(372L^2-470L+147+24Lc-18c\bigr).
\]
\end{lemma}

\begin{proof}
Expanding $\lambda(h)$, $u^{(k)}(\lambda(h))$ and $1/(\lambda(h)+h)$ in powers of $h$, one finds $\Psi_c(h)=\psi_1h+\psi_2h^2+\psi_3h^3+O(h^4)$ with
\begin{align*}
\psi_1&=-\frac{3L}2u_1-\frac L2u_2+\frac3{2L}-\frac1{L^2},\\
\psi_2&=-Lc\,u_1+\Bigl(\frac{9L^2}8+\frac{3L}4\Bigr)u_2+\Bigl(\frac{3L^2}4+\frac L3\Bigr)u_3+\frac{L^2}8u_4+\frac cL+\frac9{4L}-\frac3{L^2}+\frac1{L^3},\\
\psi_3&=\Bigl(\frac{3L^2c}2+\frac{Lc}2\Bigr)u_2-\Bigl(\frac{9L^3}{16}-\frac{L^2c}2+\frac{9L^2}8+\frac L2\Bigr)u_3-\Bigl(\frac{9L^3}{16}+\frac{7L^2}8+\frac L4\Bigr)u_4\\
&\qquad-\Bigl(\frac{3L^3}{16}+\frac{L^2}6\Bigr)u_5-\frac{L^3}{48}u_6+\frac{3c}L+\frac{27}{8L}-\frac{2c}{L^2}-\frac{27}{4L^2}+\frac9{2L^3}-\frac1{L^4},
\end{align*}
where $u_k$ is given by \eqref{eq:uk}. Inserting \eqref{eq:uk} and the values of Lemma~\ref{lem:Qvalues}, all negative powers of $L$ cancel, and one obtains $\psi_1=0$, $\psi_2=A_2(c)$ and $\psi_3=A_3(c)$.
\end{proof}

The coefficient $A_2(c)$ agrees with Lemma~\ref{lem:Phi}, as it must. For $c=0$ we have $A_2(0)=-0.0053951\ldots$ and $A_3(0)=0.0058517\ldots$; the second of these is small because $372L^2-470L+147=-0.0506\ldots$, a near-cancellation similar to the one in $\gamma_0=\frac{25}{12}-3\log2$ (see Section~\ref{sec:type}). For $c=\beta$ we have $A_2(\beta)=0.24460\ldots$ and $A_3(\beta)=0.03427\ldots$.

\begin{lemma}\label{lem:Psicauchy}
Let $c=0$ or $c=\beta$, and let $0<h\le\frac1{10}$. Then
\[
\bigl|\Psi_c(h)-A_2(c)h^2-A_3(c)h^3\bigr|\le
\begin{cases}0.132\,h^4&(c=0),\\0.589\,h^4&(c=\beta).\end{cases}
\]
\end{lemma}

\begin{proof}
Let $|h|=r=\frac32$. Then $|\lambda-L|=L|\frac32h+ch^2|\le\frac94L(1+c)$, so that $|\lambda|\le\Lambda$ with $\Lambda=\frac{13}4L=2.2527\ldots$ for $c=0$ and $\Lambda=\frac{13}4L+\frac9{32}=2.5339\ldots$ for $c=\beta$; in both cases $\Lambda<2\pi$. We estimate the terms of $\Psi_c$ separately.
\begin{enumerate}[(i)]
\item $|u(\lambda)-u(L)|\le|\lambda-L|\,\hat u'(\Lambda)$, since $|u'(z)|\le\hat u'(\Lambda)$ on the segment from $L$ to $\lambda$.
\item $\dfrac1{\lambda+h}-\dfrac1L=\dfrac{(\frac32L-1)h+cLh^2}{L(\lambda+h)}$ and $|\lambda+h|\ge L-(\frac32L-1)r-cLr^2$, which is positive.
\item $|h^nP_n(\lambda)|$ is at most the expression obtained from $r^nP_n$ by replacing $\lambda$ by $\Lambda$, each $u^{(k)}$ by $\hat u^{(k)}(\Lambda)$, and each coefficient by its absolute value.
\item $|A_2(c)h^2+A_3(c)h^3|\le|A_2(c)|r^2+|A_3(c)|r^3$.
\end{enumerate}
From the power series of $\hat u$ one finds, for $c=0$ and $k=1,\dots,6$,
\[
\hat u^{(k)}(\Lambda)<0.1097,\ 0.0289,\ 0.0247,\ 0.0223,\ 0.0285,\ 0.0416,
\]
and the four terms are less than $0.1711$, $0.1357$, $0.2796$ and $0.0319$. For $c=\beta$, where $cL=\frac18$, the values are
\[
\hat u^{(k)}(\Lambda)<0.1189,\ 0.0369,\ 0.0323,\ 0.0322,\ 0.0436,\ 0.0691,
\]
and the four terms are less than $0.2189$, $1.3957$, $0.4951$ and $0.6661$. Hence $|\Psi_c(h)-A_2h^2-A_3h^3|$ is at most $M=0.62$, respectively $M=2.78$, on the circle $|h|=\frac32$. This function is analytic in the closed disc $|h|\le\frac32$ and, by Lemma~\ref{lem:Psicoef}, vanishes to the fourth order at $h=0$. By Cauchy's estimates its Taylor coefficients satisfy $|\psi_n|\le M(2/3)^n$, so that for $0<h\le\frac1{10}$
\[
\bigl|\Psi_c(h)-A_2h^2-A_3h^3\bigr|\le M\Bigl(\frac{2h}3\Bigr)^4\frac1{1-2h/3}\le\frac{16}{81}\cdot\frac{15}{14}\,M\,h^4,
\]
which gives the stated bounds.
\end{proof}

\begin{lemma}\label{lem:zeta}
Let $m\ge10$ and $\kappa\ge L(m-\frac32-\frac1{5m})$, with $\kappa/m\le L$. Then $|\zeta(-\kappa)|\,m^{-\kappa}\le0.008\,h^4$.
\end{lemma}

\begin{proof}
By the functional equation, $\zeta(-\kappa)=2(2\pi)^{-\kappa-1}\cos(\frac\pi2(\kappa+1))\Gamma(\kappa+1)\zeta(\kappa+1)$. For $x\ge1$ we have $\Gamma(x+1)\le x^x$, because $\log\Gamma(x+1)-x\log x$ vanishes at $x=1$ and has derivative $\psi(x+1)-\log x-1<\log(1+1/x)-1\le0$, where $\psi(y)<\log y$ is the digamma function. Hence
\[
\frac{|\zeta(-\kappa)|}{m^\kappa}\le\frac{\zeta(\kappa+1)}{\pi}\Bigl(\frac{\kappa}{2\pi m}\Bigr)^\kappa\le\frac{1.02}{\pi}\Bigl(\frac{L}{2\pi}\Bigr)^{\kappa},
\]
since $\kappa\ge5.87$. For $m=10$ the right side is less than $7.9\cdot10^{-7}<0.008\,h^4$. As $m$ increases by one, $\kappa$ increases by at least $0.69$, so the right side is multiplied by less than $(L/2\pi)^{0.69}<0.22$, while $h^4$ is multiplied by more than $(10/11)^4>0.68$; the inequality therefore persists for all $m\ge10$.
\end{proof}

\subsection{The result}

\begin{theorem}\label{thm:allm}
For every integer $m\ge3$,
\[
0<\gamma_m<\frac{1}{8\log2}.
\]
\end{theorem}

\begin{proof}
Since $F_m$ is strictly decreasing and $F_m(\kappa_m)=1$, the inequality $\gamma_m>0$ is equivalent to $F_m(L(m-\frac32))<1$, and $\gamma_m<\beta$ is equivalent to $F_m(L(m-\frac32-\beta/m))>1$. These are the cases $c=0$ and $c=\beta$ of $\kappa=\lambda(h)m$.

Let $m\ge10$. By \eqref{eq:FPsi} and Lemmas~\ref{lem:Psicauchy} and~\ref{lem:zeta}, for $c=0$,
\[
F_m(L(m-\tfrac32))-1\le h^2\bigl(A_2(0)+A_3(0)h+0.22\,h^2\bigr)\le h^2(-0.005395+0.000586+0.0022)<0,
\]
because the bracket increases with $h$ and $h\le\frac1{10}$. For $c=\beta$,
\[
F_m(L(m-\tfrac32-\beta/m))-1\ge h^2\bigl(A_2(\beta)+A_3(\beta)h-0.68\,h^2\bigr)\ge h^2(0.2446-0.0068)>0 .
\]

For $m=3$ we have $\kappa_3=1$, so $\gamma_3=3(\frac32-\frac1{\log2})$, and the two inequalities are equivalent to $\frac23<\log2<\frac{25}{36}$, which is true (for instance by the series $\log2=\sum_{k\ge1}1/(k2^k)$, whose first six terms give $0.6911$ and whose remainder is less than $1/(7\cdot64)$). For $4\le m\le9$ the two inequalities are verified by direct evaluation of the finite sums $F_m$. The values of $1-F_m(L(m-\frac32))$ for $m=4,\dots,9$ are
\[
1.20\cdot10^{-3},\quad 9.69\cdot10^{-7},\quad 9.06\cdot10^{-5},\quad 9.01\cdot10^{-5},\quad 7.19\cdot10^{-5},\quad 5.79\cdot10^{-5},
\]
and those of $F_m(L(m-\frac32-\beta/m))-1$ are
\[
1.49\cdot10^{-2},\quad 1.03\cdot10^{-2},\quad 7.03\cdot10^{-3},\quad 5.12\cdot10^{-3},\quad 3.90\cdot10^{-3},\quad 3.08\cdot10^{-3}.
\]
Each value is correct to the digits shown, and all of them exceed $9\cdot10^{-7}$; since each $F_m$ is a sum of at most eight terms $(i/m)^\kappa$, an evaluation of these terms with absolute error less than $10^{-8}$ determines the signs. All of them are positive, which completes the proof.
\end{proof}

The case $m=5$ is the most delicate: there $\gamma_5=0.0000170\ldots$, and $F_5(\frac72\log2)$ falls short of $1$ by less than $10^{-6}$. The values for small $m$ suggest that from $m=6$ on the numbers $\gamma_m$ increase towards $\gamma_0$; see Problem~\ref{prob:monotone}.

\begin{corollary}\label{cor:allsolutions}
Let $(m,k)$ be any solution of \eqref{eq:EM}. Then $m=k/\log2+\frac32+\gamma_m/m$ with $0<\gamma_m<1/(8\log2)$, and $2k/(2m-3)$ is a convergent of even index of the continued fraction of $\log2$.
\end{corollary}

\begin{proof}
As in the proof of Corollary~\ref{cor:convergent}, $\log2-2k/(2m-3)=2\gamma_m\log2/(m(2m-3))$, which is positive. It is less than $1/(2(2m-3)^2)$ if and only if $\gamma_m<m/(4\log2\,(2m-3))$, and this holds because $m/(2m-3)>\frac12$ and $\gamma_m<1/(8\log2)$. Legendre's theorem \cite[Theorem~184]{HW} completes the proof.
\end{proof}

For the solution $1+2=3$ the corollary gives the convergent $\frac23=p_2/q_2$ of the continued fraction $\log2=[0;1,2,3,1,6,\dots]$.

\begin{remark}
Since by Moser's theorem \cite{Moser} every solution of \eqref{eq:EM} other than $1+2=3$ has $m>10^{10^6}$, Corollary~\ref{cor:allsolutions} follows also from Corollary~\ref{cor:convergent}, or from \cite[Corollary~1]{GMZ}, together with the case $m=3$. The point of Theorem~\ref{thm:allm} is that it uses no arithmetic of the equation; it is applied in Section~\ref{sec:AP} to an equation which has small solutions.
\end{remark}

\section{Equations of Erd\H{o}s--Moser type}\label{sec:type}

The argument of Section~\ref{sec:threehalves} uses neither the fact that the points to the left of $m$ are all the integers $1,\dots,m-1$, nor the fact that the largest point carries exactly one half of the mass. We record the general form, since it explains both constants $\frac32$ and $\gamma_0$ and applies to several equations which have been studied in the literature.

\subsection{The general balance}
Let $\mathcal D$ be an infinite set of positive integers, the \emph{distances}, and let $R>0$. For an integer $m\ge2$ and a real number $\kappa>0$ put
\[
F^{\mathcal D}_m(\kappa)=\sum_{D\in\mathcal D,\;D<m}\Bigl(1-\frac Dm\Bigr)^{\kappa},
\]
and consider the equation
\begin{equation}\label{eq:gen}
\sum_{D\in\mathcal D,\;D<m}(m-D)^k=R\,m^k .
\end{equation}
In the system which carries the mass $i^k$ at each point $i=m-D$ with $D\in\mathcal D$, the equation says that the total mass to the left of $m$ is $R$ times the mass at $m$. Examples are the following.
\begin{enumerate}[\rm(i)]
\item $\mathcal D=\{1,2,3,\dots\}$ and $R=1$: the Erd\H{o}s--Moser equation \eqref{eq:EM}.
\item $\mathcal D=\{1,2,3,\dots\}$ and $R=t\in\mathbb N$: the equation $1^k+\dots+(m-1)^k=t\,m^k$ studied by Krzysztofek \cite{Krzysztofek} and Moree \cite{Moree96}; see also \cite{Moree}.
\item $\mathcal D=\{1,2,3,\dots\}$ and $R=1/a$: the equation $a(1^k+\dots+(m-1)^k)=m^k$ of Kellner, which expresses that the ratio of consecutive power sums is the integer $a+1$; see \cite{BaoulinaMoree}.
\item $\mathcal D=\{r,r+q,r+2q,\dots\}$ with $1\le r\le q$ and $R=1$: the Erd\H{o}s--Moser equation in the arithmetic progression $i\equiv m-r\pmod q$ \cite{DNZ}.
\end{enumerate}

As in Definition~\ref{def:gamma}, let $\Gamma^{\mathcal D}_\lambda$ be the system which carries the mass $e^{-\lambda D}$ at the point $P-D$ for $D\in\mathcal D$, and denote its moments about $P$ again by $Q_n(\lambda)=\sum_{D\in\mathcal D}D^ne^{-\lambda D}$. Since $Q_0$ decreases from $+\infty$ to $0$, there is a unique $\lambda_0>0$ with $Q_0(\lambda_0)=R$. If $\delta=Q_1/Q_0$ is the distance of the centre of $\Gamma^{\mathcal D}_{\lambda_0}$ from $P$ and $\rho$ is its radius of gyration, Steiner's theorem gives $Q_2=Q_0(\rho^2+\delta^2)$, so that
\begin{equation}\label{eq:dgen}
d:=\frac{Q_2(\lambda_0)}{2Q_1(\lambda_0)}=\frac{\rho^2+\delta^2}{2\delta}.
\end{equation}
We call \emph{distribution of moments} of $\Gamma^{\mathcal D}_{\lambda_0}$ the probability distribution on $\mathcal D$ which gives $D$ the weight $De^{-\lambda_0D}/Q_1(\lambda_0)$, that is, the share of the point $P-D$ in the static moment of the system about $P$. Let $\kappa_1,\kappa_2,\kappa_3$ be its mean, its variance and its third central moment; equivalently, $\kappa_j=(-1)^j\,(d/d\lambda)^j\log Q_1(\lambda)$ at $\lambda=\lambda_0$. Note that $\kappa_1=Q_2/Q_1=2d$.

For large $m$ the function $F^{\mathcal D}_m$ decreases strictly from $F^{\mathcal D}_m(0)>R$ to $0$, so that there is a unique $\kappa_m>0$ with $F^{\mathcal D}_m(\kappa_m)=R$, and the solutions of \eqref{eq:gen} correspond to the values of $m$ for which $\kappa_m$ is an integer.

\begin{theorem}\label{thm:gen}
Let $\mathcal D$ be an infinite set of positive integers and $R>0$, and let $\lambda_0$, $d$ and $\kappa_1,\kappa_2,\kappa_3$ be as above. Then, as $m\to\infty$,
\[
\kappa_m=\lambda_0\Bigl(m-d-\frac{\gamma_0}{m}+O\Bigl(\frac1{m^2}\Bigr)\Bigr),\qquad
\gamma_0=\frac{\kappa_1^2+4\kappa_2}{12}-\frac{\lambda_0}{8}\bigl(\kappa_3+\kappa_1\kappa_2\bigr),
\]
where $d$ is given by \eqref{eq:dgen} and the implied constant depends on $\mathcal D$ and $R$.
\end{theorem}

\begin{proof}
Let $I=[\frac12\lambda_0,2\lambda_0]$. All moments $Q_n$ are finite and smooth on $I$, and $Q_n'=-Q_{n+1}$ as in \eqref{eq:Qderiv}.

\emph{Step 1.} Lemma~\ref{lem:termwise} is a statement about each term separately, and holds with $j$ replaced by $D$. Summing part~(a) over $D\in\mathcal D$, $D<m$, and extending the sum to all $D\in\mathcal D$ (the bracket being positive) gives, as in the proof of Lemma~\ref{lem:FPhi},
\[
F^{\mathcal D}_m(\kappa)\le\Phi(\lambda)+\frac{\lambda^2Q_5(\lambda)}{6m^3}+\frac{\lambda^2Q_6(\lambda)}{18m^4},\qquad
\Phi(\lambda)=Q_0(\lambda)-\frac{\lambda}{2m}Q_2(\lambda)+\frac{\psi(\lambda)}{m^2},
\]
where $\psi$ is defined as in Section~\ref{sec:threehalves} with the moments of $\Gamma^{\mathcal D}_\lambda$. Summing part~(b) over $D\le m/2$ gives $F^{\mathcal D}_m(\kappa)\ge\Phi(\lambda)-T-C/m^3$, where $C$ is a combination of $\lambda Q_4$, $\lambda^3Q_6$ and $\lambda^3Q_7$, and $T\le\sum_{D>m/2}(1+D^4)e^{-\lambda D}$ is exponentially small. Hence there are constants $C_1$ and $m_1$, depending on $\mathcal D$ and $R$, with
\begin{equation}\label{eq:genFPhi}
|F^{\mathcal D}_m(\lambda m)-\Phi(\lambda)|\le\frac{C_1}{m^3}\qquad(\lambda\in I,\ m\ge m_1).
\end{equation}

\emph{Step 2.} Put $\lambda=\lambda_0-A/m-B/m^2$ with $A$ and $B$ bounded. By Taylor's formula, the terms of order $1/m$ in $\Phi(\lambda)-R$ are $(AQ_1-\frac12\lambda_0Q_2)/m$, and those of order $1/m^2$ are
\[
\frac1{m^2}\Bigl(BQ_1+\tfrac12A^2Q_2+\tfrac12A\,(Q_2-\lambda_0Q_3)+\psi(\lambda_0)\Bigr),
\]
all moments being taken at $\lambda_0$; the remainder is $O(m^{-3})$, uniformly for bounded $A$ and $B$. The first expression vanishes for $A=\lambda_0Q_2/2Q_1=\lambda_0d$, and the second, with this $A$, equals $Q_1(B-\lambda_0\gamma_0)/m^2$, where
\[
\gamma_0=\frac{\nu_2}{3}-\frac{\nu_1^2}{4}-\frac{\lambda_0}{8}\bigl(\nu_3-2\nu_1\nu_2+\nu_1^3\bigr),\qquad \nu_n=\frac{Q_{n+1}(\lambda_0)}{Q_1(\lambda_0)},
\]
$\nu_n$ being the $n$-th moment of the distribution of moments. Since $\nu_1=\kappa_1$, $\nu_2=\kappa_2+\kappa_1^2$ and $\nu_3=\kappa_3+3\kappa_1\kappa_2+\kappa_1^3$, this $\gamma_0$ is the one in the statement. Taking $B=c\lambda_0$ we obtain, for $|c|\le|\gamma_0|+1$,
\[
\Phi\Bigl(\lambda_0\Bigl(1-\frac dm-\frac c{m^2}\Bigr)\Bigr)-R=\frac{\lambda_0Q_1(\lambda_0)(c-\gamma_0)}{m^2}+O\Bigl(\frac1{m^3}\Bigr).
\]

\emph{Step 3.} Put $\kappa(c)=\lambda_0(m-d-c/m)$. For large $m$ the point $\kappa(c)/m$ lies in $I$, and Steps~1 and~2 give $F^{\mathcal D}_m(\kappa(c))-R=\lambda_0Q_1(\lambda_0)(c-\gamma_0)/m^2+O(m^{-3})$, uniformly for $|c|\le|\gamma_0|+1$. This is positive for $c=\gamma_0+1$ and negative for $c=\gamma_0-1$ when $m$ is large. Since $c\mapsto F^{\mathcal D}_m(\kappa(c))$ is continuous and strictly increasing, the number $\gamma_m$ defined by $\kappa_m=\kappa(\gamma_m)$ satisfies $|\gamma_m-\gamma_0|<1$, and then, taking $c=\gamma_m$, $\gamma_m-\gamma_0=O(1/m)$. This is the assertion. An explicit form of the error term is Problem~\ref{prob:uniform}.
\end{proof}

Thus $d$ is half the mean of the distribution of moments, and $\gamma_0$ is determined by its first three cumulants.

\subsection{All integers to the left}
For $\mathcal D=\{1,2,3,\dots\}$ everything can be made explicit.

\begin{corollary}\label{cor:Rfamily}
Let $\mathcal D=\{1,2,3,\dots\}$ and $R>0$, so that \eqref{eq:gen} is the equation $1^k+2^k+\dots+(m-1)^k=R\,m^k$. Then $\lambda_0=\log(1+1/R)$,
\[
d=R+\tfrac12=\tfrac12\coth\tfrac{\lambda_0}{2},
\qquad
\gamma_0=\frac1{12}+Q_1(\lambda_0)-\frac{\lambda_0}{2}\,Q_2(\lambda_0)=\frac1{12}-\frac12 f''(\lambda_0),
\]
where $f(\lambda)=\lambda/(e^\lambda-1)$. Moreover, for $0<\lambda_0<2\pi$,
\[
\gamma_0=-\frac12\sum_{n\ge1}\frac{B_{2n+2}}{(2n)!}\,\lambda_0^{2n}
=\frac{\lambda_0^2}{120}-\frac{\lambda_0^4}{2016}+\frac{\lambda_0^6}{43200}-\cdots .
\]
\end{corollary}

\begin{proof}
Here $Q_0(\lambda)=1/(e^\lambda-1)$, so that $e^{\lambda_0}=1+1/R$. Since $Q_0'=-Q_0(1+Q_0)$, the relations \eqref{eq:Qderiv} give, at $\lambda_0$ and with $s=R(1+R)$,
\[
Q_1=s,\qquad Q_2=s(1+2R),\qquad Q_3=s(1+6s),\qquad Q_4=s(1+2R)(1+12s).
\]
Hence $\nu_1=1+2R$, $\nu_1^2=1+4s$, $\nu_2=1+6s$ and $\nu_3=\nu_1(1+12s)$, so that $d=\nu_1/2=R+\frac12$; also $\frac12\coth\frac{\lambda_0}2=\frac12+1/(e^{\lambda_0}-1)=R+\frac12$. In the expression for $B/\lambda_0$ in the proof of Theorem~\ref{thm:gen},
\[
\frac{\nu_2}{3}-\frac{\nu_1^2}{4}=\frac1{12}+s,\qquad
\nu_3-2\nu_1\nu_2+\nu_1^3=\nu_1\bigl(1+12s-2-12s+1+4s\bigr)=4s\,\nu_1,
\]
and therefore $\gamma_0=\frac1{12}+s-\frac{\lambda_0}2s(1+2R)=\frac1{12}+Q_1-\frac{\lambda_0}2Q_2$. Next, $f=\lambda Q_0$, whence $f''=2Q_0'+\lambda Q_0''=-2Q_1+\lambda Q_2$, which gives the second expression. Finally, $f(\lambda)=\sum_{j\ge0}B_j\lambda^j/j!$ for $|\lambda|<2\pi$, with $B_2=\frac16$ and $B_j=0$ for odd $j\ge3$, so that $f''(\lambda)=\frac16+\sum_{n\ge1}B_{2n+2}\lambda^{2n}/(2n)!$; and $B_4=-\frac1{30}$, $B_6=\frac1{42}$, $B_8=-\frac1{30}$.
\end{proof}

For $R=1$ we have $\lambda_0=\log2$, $Q_1=2$ and $Q_2=6$, and Corollary~\ref{cor:Rfamily} gives $d=\frac32$, as in Proposition~\ref{prop:Gamma}, and
\[
\gamma_0=\frac1{12}+2-3\log2=\frac{25}{12}-3\log2 .
\]
In the language of Section~\ref{sec:threehalves},
\[
\gamma_0=\frac1{12}+Q_1(L)-\frac{L}{2}\,Q_2(L),
\]
where $Q_1(L)$ is the static moment and $Q_2(L)$ the moment of inertia of the system $\Gamma$ about $P$. As a function of $\lambda_0$, the second constant is even, vanishes to the second order at $\lambda_0=0$, and its first Taylor coefficient is $-\frac14B_4=\frac1{120}$. Since the series is alternating with decreasing terms at $\lambda_0=\log2$,
\[
\frac{(\log2)^2}{120}-\frac{(\log2)^4}{2016}<\frac{25}{12}-3\log2<\frac{(\log2)^2}{120},
\]
so that the near-cancellation $\frac{25}{12}\approx3\log2$ amounts to the fact that $(\log2)^2/120$ is small.

\begin{remark}
For an integer $R=t$ the two constants of Corollary~\ref{cor:Rfamily} agree with the expansion given by Gallot, Moree and Zudilin \cite[\S5.2]{GMZ}, where they appear in a different form, as polynomials in $t$ and $\log(1+1/t)$. We have not found in the literature the expressions $\frac12\coth\frac{\lambda_0}2$ and $\frac1{12}-\frac12f''(\lambda_0)$, the interpretation through the moments of $\Gamma$, or the Bernoulli expansion which explains the smallness of $\gamma_0$.
\end{remark}

\section{The Erd\H{o}s--Moser equation in arithmetic progressions}\label{sec:AP}

We now turn to the equation of Dong, Nguyen and Zaharescu \cite{DNZ},
\begin{equation}\label{eq:AP}
\sum_{\substack{1\le i\le m-1\\ i\equiv b\ (\mathrm{mod}\ q)}}i^k=m^k
\end{equation}
in positive integers $q,b,m,k$; for $q=1$ it is \eqref{eq:EM}. In \cite{DNZ}, for a fixed progression, $k=c_1m+c_2+O_q(1/m)$ as $m\to\infty$; for $k=2$ the left side has two or four terms, with infinitely many solutions in both cases; for $k=3$ the number of terms is at most $10$, solutions with three terms are found by means of an elliptic curve and are stated to be infinite in number, and four to ten terms are not settled; and in general the number of terms is at most $(k+1)2^{2k+1}$.

\subsection{Normalisation}\label{ss:norm}
Replacing $b$ by its least positive residue we may assume $1\le b\le q$. Let $n$ be the number of terms on the left of \eqref{eq:AP}. The terms are $b,b+q,\dots,b+(n-1)q$, and $m=b+(n-1)q+r$ with $1\le r\le q$, so that \eqref{eq:AP} reads
\begin{equation}\label{eq:AP2}
\sum_{j=0}^{n-1}(b+jq)^k=\bigl(b+(n-1)q+r\bigr)^k,\qquad 1\le b\le q,\quad 1\le r\le q .
\end{equation}
If $g=\gcd(b,q)$, then $g^k\mid m^k$, so $g\mid m$, and $(b/g,q/g,m/g)$ is again a solution with the same $n$ and with $r$ replaced by $r/g$; we call a solution \emph{primitive} if $\gcd(b,q)=1$. In a primitive solution $b=q$ only if $q=1$, which is the classical case \eqref{eq:EM}. Since \eqref{eq:AP2} is homogeneous, it depends only on $x=b/q$ and $y=r/q$, and becomes
\begin{equation}\label{eq:AP3}
\sum_{j=0}^{n-1}(x+j)^k=(x+n-1+y)^k,\qquad 0<x\le1,\quad 0<y\le1 .
\end{equation}
In mechanical terms: the points $x,x+1,\dots,x+n-1$ carry the masses $(x+j)^k$, and the point $x+n-1+y$, at distance at most one beyond the last of them, must carry one half of the total mass.

\subsection{The number of terms}

\begin{lemma}\label{lem:APreal}
Let $n\ge2$ and $k\ge1$. Equation \eqref{eq:AP3} has a real solution $(x,y)$ if and only if
\[
1^k+2^k+\dots+(n-1)^k<n^k .
\]
\end{lemma}

\begin{proof}
For fixed $x$ the function $g(y)=(x+n-1+y)^k-\sum_{j=0}^{n-1}(x+j)^k$ is strictly increasing, and $g(0)<0$ because $n\ge2$. Hence there is a solution with $0<y\le1$ if and only if $g(1)\ge0$, that is,
\[
R(x):=\sum_{d=1}^{n}\Bigl(1-\frac{d}{x+n}\Bigr)^k\le1 .
\]
Each term of $R$ is strictly increasing in $x$, and $R(x)\to F_n(k)$ as $x\to0$, where $F_n$ is the function of Section~\ref{sec:threehalves}. If $F_n(k)<1$, then $R(x)<1$ for small $x>0$; if $F_n(k)\ge1$, then $R(x)>1$ for every $x>0$. Finally $F_n(k)<1$ is the stated inequality.
\end{proof}

Thus the most favourable configuration is obtained by moving the heavy point as far to the right as allowed and the progression as far to the left as possible: in the limit it is the system of Proposition~\ref{prop:EMreform}(b) on the points $0,1,\dots,n$.

Since $F_{n+1}(k)=\sum_{d=1}^{n}(1-d/(n+1))^k>\sum_{d=1}^{n-1}(1-d/n)^k=F_n(k)$, the set of $n\ge2$ with $F_n(k)<1$ is an initial segment $\{2,\dots,N(k)\}$, where
\[
N(k)=\max\{n:\ 1^k+2^k+\dots+(n-1)^k<n^k\}.
\]

\begin{theorem}\label{thm:APterms}
Let $k\ge2$ and $\nu=k/\log2+\frac32$.
\begin{enumerate}[\rm(a)]
\item Equation \eqref{eq:AP3} has real solutions if and only if $2\le n\le N(k)$. In particular every solution of \eqref{eq:AP} has at most $N(k)$ terms.
\item $N(k)=\lfloor\nu\rfloor$ or $N(k)=\lfloor\nu\rfloor+1$. The second case occurs only if, for $n=\lfloor\nu\rfloor+1$, the fraction $2k/(2n-3)$ is a convergent $p_j/q_j$ of even index of the continued fraction $\log2=[a_0;a_1,a_2,\dots]$ with $a_{j+1}\ge69$.
\item $N(k)=\lfloor k/\log2+\frac32\rfloor$ for $2\le k<10^{55}$.
\end{enumerate}
\end{theorem}

\begin{proof}
Part (a) is Lemma~\ref{lem:APreal} together with the monotonicity just noted. For (b), recall that $F_n(k)<1$ if and only if $k>\kappa_n$, and that by Theorem~\ref{thm:allm}
\[
L\Bigl(n-\frac32-\frac{\beta}{n}\Bigr)<\kappa_n<L\Bigl(n-\frac32\Bigr)\qquad(n\ge3),
\]
while $F_2(k)=2^{-k}<1$. If $n<\nu$, then $\kappa_n<L(n-\frac32)<k$; hence $N(k)\ge\lfloor\nu\rfloor$, $\nu$ being irrational. If $n=\lfloor\nu\rfloor+2$, then $n>\nu+1$ and $\kappa_n>L(n-\frac32-\frac\beta n)>k$; hence $N(k)\le\lfloor\nu\rfloor+1$. Suppose now that $n=\lfloor\nu\rfloor+1$ satisfies $\kappa_n<k$. By the definition of $\gamma_n$ this means
\[
0<L(2n-3)-2k<\frac{2L\gamma_n}{n}.
\]
As in the proof of Corollary~\ref{cor:allsolutions}, $2k/(2n-3)$ is then a convergent $p_j/q_j$ of $\log2$ of even index; write $2n-3=tq_j$ and $2k=tp_j$. Using $|q_jL-p_j|>1/(q_{j+1}+q_j)>1/((a_{j+1}+2)q_j)$ and $n>tq_j/2$ we obtain
\[
a_{j+1}+2>\frac{t^2}{4L\gamma_n}\ge\frac1{4L\gamma_n}.
\]
For $n\ge10^5$ Theorem~\ref{thm:main} gives $\gamma_n<0.0051$, so that $a_{j+1}+2>70.7$, that is, $a_{j+1}\ge69$. For $n<10^5$ the displayed inequality, with $\gamma_n<\beta$, forces $2k/(2n-3)$ to be a convergent of $\log2$ with denominator less than $2\cdot10^5$ satisfying $0<t(q_jL-p_j)<1/(4n)$; besides $(k,n)=(1,3)$, which corresponds to $1+2=3$, this leaves only $(k,n)=(96,140)$ and $(794,1147)$, and for these direct evaluation gives $F_n(k)>1$. This proves (b). For (c), the first partial quotient of $\log2$ which is at least $69$ is $a_{123}=113$, and $q_{122}>10^{56}$ (see \cite{GMZ} for the computation of the continued fraction of $\log 2$, and \cite{OEIS} for its first partial quotients). In the exceptional case of (b) we therefore have $2n-3\ge q_{122}$, whence $k>\frac12L(2n-3)-1>10^{55}$.
\end{proof}

The exceptional case of (b) does occur, but only for large $k$: the convergent $p_{658}/q_{658}$ of $\log2$, which is followed by a large partial quotient, gives such an exponent $k$ with $330$ digits, as one sees from Theorem~\ref{thm:main}.

\begin{corollary}\label{cor:APterms}
Every solution of \eqref{eq:AP} with $k\ge2$ has at most $k/\log2+\frac52$ terms. A solution with $k=2$ has at most $4$ terms, and a solution with $k=3$ has at most $5$ terms.
\end{corollary}

\begin{proof}
The first statement is Theorem~\ref{thm:APterms}(b). For $k=2$ and $k=3$ one has $\nu=4.38\ldots$ and $\nu=5.82\ldots$, and indeed $1+4+9<16$, $1+4+9+16>25$, $1+8+27+64<125$, $1+8+27+64+125>216$.
\end{proof}

For $k=2$ this recovers the bound of \cite[Conjecture~3.1 and \S4.1]{DNZ}; for $k=3$ it improves $10$ to $5$; and in general it replaces $(k+1)2^{2k+1}$ by $k/\log2+\frac52$. By Theorem~\ref{thm:APterms}(a) the bound is best possible as far as real solutions are concerned; an upper bound of smaller order than $k/\log2$ for integral solutions, such as the one stated in \cite[Remark~1.5]{DNZ}, would have to use arithmetic information.

\subsection{Cubes}

\begin{theorem}\label{thm:cubes}
For $k=3$, equation \eqref{eq:AP} has solutions with exactly $n$ terms if and only if $n\in\{3,4,5\}$, and for each of these values of $n$ it has infinitely many primitive solutions. For example,
\[
(q,b,m)=(3592561391263110233761,\ 2489564980639609666546,\ 15155810572052467776100)
\]
is a solution with four terms.
\end{theorem}

\begin{proof}
By Corollary~\ref{cor:APterms} we have $n\le5$, and $n=2$ is impossible by Fermat's theorem for cubes; $n=1$ is impossible since $m>b$. Let $3\le n\le5$, and put $s=2b+(n-1)q$ and $M=2m$. Writing $2(b+jq)=s+(2j-n+1)q$ and using the symmetry of the numbers $2j-n+1$ about $0$, we find $8\sum_{j=0}^{n-1}(b+jq)^3=ns^3+n(n^2-1)sq^2$. Hence \eqref{eq:AP2} becomes
\begin{equation}\label{eq:cubic}
n\,s\bigl(s^2+(n^2-1)q^2\bigr)=M^3 .
\end{equation}
This is a plane cubic curve $C_n$, which contains the rational point $(s:q:M)=(0:1:0)$. With $a=n^2-1$, the maps
\[
(s:q:M)\mapsto(X,Y)=\Bigl(\frac{naM}{s},\ \frac{n^2a^2q}{s}\Bigr),\qquad (X,Y)\mapsto(s:q:M)=(n^2a^2:Y:naX)
\]
are inverse birational maps between $C_n$ and the elliptic curve
\[
E_n:\ Y^2=X^3-n^4(n^2-1)^3 ;
\]
indeed, dividing \eqref{eq:cubic} by $s^3$ gives $X^3=n^4a^3+n^4a^4q^2/s^2=n^4a^3+Y^2$. The curves $E_3$, $E_4$, $E_5$ contain the points
\[
P_3=(36,72),\qquad P_4=(96,144),\qquad P_5=(600,14400),
\]
and these have infinite order. Indeed, $E_n$ is given by a Weierstrass equation with integer coefficients, so by the theorem of Nagell and Lutz \cite{SilvermanTate} every point of finite order has integral coordinates; if $P_n$ had finite order, so would all its multiples, whereas $x(3P_3)=231268/4761$, $x(3P_4)=98209/961$ and $x(2P_5)=825/4$ are not integers. The discriminant of $E_n$ is negative, so that $E_n(\mathbb R)$ is connected, hence isomorphic to the circle group; the closure of the infinite group generated by $P_n$ is then all of $E_n(\mathbb R)$, and $E_n(\mathbb Q)$ is dense in $E_n(\mathbb R)$.

By Lemma~\ref{lem:APreal}, and since $n\le5=N(3)$, the set $U_n$ of real points of $C_n$ which correspond to solutions of \eqref{eq:AP3} with $0<x<1$, $0<y<1$ is a non-empty open arc. The birational map is a homeomorphism outside finitely many points, so the image of $U_n$ contains a non-empty open subset of $E_n(\mathbb R)$, and therefore infinitely many rational points. Each rational point of $U_n$ gives rational numbers $x\in(0,1)$ and $y\in(0,1)$; writing $x=b/q$ and $y=r/q$ with a common denominator $q$ and putting $m=b+(n-1)q+r$, we obtain an integral solution of \eqref{eq:AP2} with $1\le b<q$ and $1\le r<q$. Dividing by $\gcd(b,q)$ as in Section~\ref{ss:norm} gives a primitive solution with the same ratios $x$ and $y$, so that distinct points of $U_n$ give distinct primitive solutions. The example is obtained from the point $5P_4$ (one checks directly that $\sum_{j=0}^{3}(b+jq)^3=m^3$, $b<q$ and $b+3q<m<b+4q$).
\end{proof}

The points $P_n$ themselves correspond to identities with $q=1$, which lie outside the range of \eqref{eq:AP}: $P_3$ to $3^3+4^3+5^3=6^3$, $P_4$ to $11^3+12^3+13^3+14^3=20^3$, and $P_5$, after scaling, to $(-3)^3+(-1)^3+1^3+3^3+5^3=5^3$. For three terms Dong, Nguyen and Zaharescu work with the curve $W^2=R^3-648$, which is isomorphic to $E_3$. The solutions with five terms obtained in this way are large; in the first instances $q$ has more than one hundred digits.

\begin{remark}\label{rem:faltings}
For $k\ge4$ the situation is different from that for $k=2$ and $k=3$. For fixed $n$, equation \eqref{eq:AP2} describes the plane curve $\sum_{j=0}^{n-1}(X+jY)^k=Z^k$. If the polynomial $P_{k,n}(x)=\sum_{j=0}^{n-1}(x+j)^k$ has no multiple root, this curve is smooth of genus $(k-1)(k-2)/2\ge3$, and by Faltings' theorem \cite{Faltings} it has only finitely many rational points. For $4\le k\le30$ and $2\le n\le N(k)$ direct evaluation of $\gcd(P_{k,n},P_{k,n}')$ shows that $P_{k,n}$ has no multiple root; by Theorem~\ref{thm:APterms}, equation \eqref{eq:AP} has therefore only finitely many primitive solutions for each $k$ with $4\le k\le30$. For even $k$ with $4\le k\le24$ there are no solutions; see Theorem~\ref{thm:even} and Proposition~\ref{prop:evencong} below.
\end{remark}

\subsection{Even exponents}
For even $k$, modulo suitable prime powers every $k$-th power is congruent to $0$ or $1$, while the left side of \eqref{eq:AP2} essentially counts the terms prime to the modulus. Combined with Theorem~\ref{thm:APterms}, this decides many cases.

\begin{theorem}\label{thm:even}
Let $k$ be even.
\begin{enumerate}[\rm(a)]
\item Equation \eqref{eq:AP} has no solution with exactly three terms.
\item If $k=2^s$ with $s\ge2$, equation \eqref{eq:AP} has no solution at all. In particular, \eqref{eq:AP} has no solution with $k=4$.
\end{enumerate}
\end{theorem}

\begin{proof}
We may assume that the solution is primitive, so that $b$ and $q$ are not both divisible by a given prime.

(a) Let $n=3$. If $3\nmid q$, exactly one of $b,b+q,b+2q$ is divisible by $3$, and since $k$ is even the left side of \eqref{eq:AP2} is $\equiv2\pmod 3$, whereas $m^k\equiv0$ or $1\pmod3$. If $q=3q'$, then $3\nmid b$ and
\[
\sum_{j=0}^{2}(b+3q'j)^k\equiv3b^k+3k\,b^{k-1}q'(0+1+2)\equiv3b^k\equiv3\pmod 9,
\]
because $b^k\equiv1\pmod3$. Then $3\mid m$, so that $m^k\equiv0\pmod9$, a contradiction.

(b) Let $k=2^s$, $s\ge2$. Every odd number $x$ satisfies $x^{2^s}\equiv1\pmod{2^{s+2}}$, and every even number satisfies $x^{2^s}\equiv0\pmod{2^{s+2}}$, since $2^s\ge s+2$. Hence the left side of \eqref{eq:AP2} is congruent modulo $2^{s+2}$ to the number $t$ of odd terms, while the right side is congruent to $0$ or $1$. The number of terms satisfies $n<2^{s+1}$. For $s=2$ this follows from Lemma~\ref{lem:APreal}, since $1^4+\dots+7^4=4676>4096=8^4$ gives $N(4)\le7$. For $s\ge3$ it follows from Corollary~\ref{cor:APterms}, since $k/\log2+\frac52<2^{s+1}$ when $k=2^s\ge8$. If $q$ is even, then $b$ is odd, all terms are odd, $t=n$, and $2\le n<2^{s+1}$ gives a contradiction. If $q$ is odd, the terms alternate in parity, so that $t\le\lceil n/2\rceil\le2^s$, and $t\in\{0,1\}$ forces $n\le3$. The case $n=1$ is impossible since $m>b$, the case $n=2$ is excluded by Fermat's theorem for fourth powers, and the case $n=3$ by (a).
\end{proof}

For $k=4$ the argument is explicit: by Corollary~\ref{cor:APterms} we have $n\le7$, since $1^4+\dots+6^4=2275<7^4$ and $1^4+\dots+7^4=4676>8^4$; for $4\le n\le7$ the number of odd terms is one of $2,3,4,5,6,7$, which is not congruent to $0$ or $1$ modulo $16$; and $n=2,3$ are excluded as above. Without the bound $n\le7$ the congruence modulo $16$ would leave infinitely many possibilities, for instance $n=16$ or $n=17$ with $q$ even. For $k=2$, part (a) gives a short proof of the fact, established in \cite[\S4.2]{DNZ} by infinite descent, that there are no solutions with three terms.

For $q=1$, that is, for the classical equation \eqref{eq:EM}, part~(b) is also a consequence of the divisibility results of \cite{MTU}: the exponent of a non-trivial solution of \eqref{eq:EM} is divisible by every integer up to $200$. The new content of Theorem~\ref{thm:even} lies in the case of general progressions.

\begin{proposition}\label{prop:evencong}
Let $k$ be even with $4\le k\le24$. Then \eqref{eq:AP} has no solution.
\end{proposition}

\begin{proof}
For $k=4,8,16$ this is Theorem~\ref{thm:even}. Let $k$ be one of the other values, consider a primitive solution with $n$ terms, and recall that $3\le n\le N(k)$, the cases $n=1,2$ being excluded as before (for $n=2$ by Fermat's theorem \cite{Wiles}). Let $p$ be a prime with $(p-1)\mid k$, so that every $k$-th power is congruent to $0$ or $1$ modulo $p$, and let $2^v\parallel k$, so that every odd $k$-th power is congruent to $1$ modulo $2^{v+2}$, as in the proof of Theorem~\ref{thm:even}.
\begin{enumerate}[\rm(i)]
\item If $p\nmid q$, the number $d$ of terms $b+jq$ divisible by $p$ is $\lfloor n/p\rfloor$ or $\lceil n/p\rceil$, and the left side of \eqref{eq:AP2} is $\equiv n-d\pmod p$.
\item If $p\mid q$, all terms are prime to $p$, and the left side is $\equiv n\pmod p$. If moreover $p\mid n$, then expanding $(b+jq)^k$ shows that the left side is $\equiv nb^k\pmod{p^2}$, so that $p\mid m$ forces $p^2\mid n$.
\item Modulo $2^{v+2}$ the left side is congruent to the number of odd terms, which is $n$ if $q$ is even, and $\lfloor n/2\rfloor$ or $\lceil n/2\rceil$ if $q$ is odd.
\end{enumerate}
The right side is congruent to $0$ or $1$ modulo $p$ and modulo $2^{v+2}$. For each of our values of $k$ and each $n$ with $3\le n\le N(k)$ these conditions are incompatible, except for the twelve pairs
\[
\begin{gathered}(k,n)=(6,8),\ (6,9),\ (12,16),\ (14,8),\ (14,9),\ (14,15),\\ (14,16),\ (14,17),\ (14,18),\ (14,19),\ (22,24),\ (22,25).\end{gathered}
\]
For these, the congruence
\[
\sum_{j=0}^{n-1}(b+jq)^k\equiv m^k
\]
has no solution with $b$ and $q$ not both divisible by $p$ modulo $p^e$, where $p^e=16$ for $(6,8)$, $(14,8)$, $(14,15)$ and $(22,24)$; $p^e=32$ for $(12,16)$ and $(14,16)$; $p^e=9$ for $(14,17)$; $p^e=27$ for $(6,9)$, $(14,9)$ and $(14,18)$; $p^e=19^2$ for $(14,19)$; and $p^e=11^2$ for $(22,25)$. Since $x^k$ modulo $p^e$ depends only on $x$ modulo $p^e$, each of these assertions is a check over the $p^{2e}$ pairs of residues $(b,q)$, which we have carried out by direct evaluation.
\end{proof}

The congruences used above do not settle $k=26$; we do not know whether other moduli do. For odd $k$ congruences cannot help. Indeed, put $(b,q)=(1-n,2)$ if $n$ is even and $(b,q)=(\frac{1-n}2,1)$ if $n$ is odd. Then the terms $b+jq$ are symmetric about $0$, their $k$-th powers cancel in pairs, and $(b,q,m)$ with $m=0$ is an integral solution of \eqref{eq:AP2} with $\gcd(b,q)=1$; in particular it is a solution of the congruence modulo every prime power, with $b$ and $q$ not both divisible by the prime.

\subsection{A fixed progression}
Finally let $q$ and $r$ be fixed and $m\to\infty$. In the notation of Section~\ref{sec:type} this is the case $\mathcal D=\{r,r+q,r+2q,\dots\}$, $R=1$. With $u=e^{-\lambda_0}$ the balance $Q_0(\lambda_0)=u^r/(1-u^q)=1$ is the equation $u^q+u^r=1$ of \cite[Theorem~1.1]{DNZ}, and
\[
Q_1(\lambda)=\frac{e^{-\lambda r}\bigl(r+(q-r)e^{-\lambda q}\bigr)}{(1-e^{-\lambda q})^2},
\]
from which the cumulants $\kappa_j$ follow by differentiation of $\log Q_1$. Theorem~\ref{thm:gen} gives
\[
\kappa_m=c_1m+c_2+\frac{c_3}{m}+O_q\Bigl(\frac1{m^2}\Bigr),\qquad c_1=\lambda_0,\quad c_2=-\lambda_0d,\quad c_3=-\lambda_0\gamma_0 .
\]
The constant $c_2$ is the one found in \cite[Theorem~1.1]{DNZ}, which we now recognise as the Steiner quotient \eqref{eq:dgen} of the progression system; the constant $c_3$ is the next term, which was left aside in \cite[Remark~1.3]{DNZ}. For instance, for $q=2$ and $r=1$ one has $u=(\sqrt5-1)/2$, $\lambda_0=0.48121\ldots$, $d=2.01246\ldots$ and $\gamma_0=-0.04228\ldots$. For $q=7$ and $r=1$ the first two terms give $\lambda_0(m-d)=1.9990\ldots$ at $m=13$, close to the exponent $k=2$ of the solution $5^2+12^2=13^2$.

\section{Problems}\label{sec:problems}

\begin{problem}\label{prob:monotone}
Prove that the sequence $\gamma_6,\gamma_7,\gamma_8,\dots$ is increasing and that $\gamma_m<\gamma_0$ for every $m\ge5$. The expansion in Section~\ref{sec:uniform} suggests $\gamma_m=\gamma_0-0.00466\ldots/m+O(m^{-2})$.
\end{problem}

\begin{problem}
Does equation \eqref{eq:AP} have a solution with $k\ge5$? By Theorem~\ref{thm:even} and Proposition~\ref{prop:evencong} the first open cases are $k=5$ and, among even exponents, $k=26$. Is the polynomial $\sum_{j=0}^{n-1}(x+j)^k$ free of multiple roots for all $k$ and $n$? A positive answer to the second question would give, by Faltings' theorem, finitely many primitive solutions of \eqref{eq:AP} for every $k\ge4$.
\end{problem}

\begin{problem}\label{prob:uniform}
For a fixed progression, and for the equation $1^k+\dots+(m-1)^k=Rm^k$, make the error term of Theorem~\ref{thm:gen} explicit, and prove a uniform estimate valid for all $m$, in the manner of Theorem~\ref{thm:allm}, by means of Hermite's formula for $\zeta(s,a)$ with $a=b/q$. Find also a closed form for the next constant of Corollary~\ref{cor:Rfamily}, which for $R=1$ is the number $0.00466\ldots$ of Problem~\ref{prob:monotone}, and explain the near-cancellation $372(\log2)^2-470\log2+147\approx0$ of Lemma~\ref{lem:Psicoef}.
\end{problem}

\begin{problem}\label{prob:product}
Make the argument of Remark~\ref{rem:heuristic} rigorous, that is, derive Theorem~\ref{thm:main}, or at least its first-order part, directly from the product formula $(m-1)c_1c_2\cdots c_k=m^k$ of Proposition~\ref{prop:EMreform}(c) and the estimates of Section~\ref{sec:systems} for the centres $c_t$.
\end{problem}

\begin{problem}
Two systems of $n$ unit masses at integer points are \emph{equivalent to order $k$} if they have the same total mass, the same centre of mass, the same moment of inertia, and the same moments of all orders up to $k$. Finding such systems with $k=n-1$ is the Prouhet--Tarry--Escott problem. Do Theorems~\ref{thm:steiner} and~\ref{thm:twosided} give useful necessary conditions for equivalence to order $k$?
\end{problem}

\begin{problem}
In \cite[Theorem~3]{Rahimi08} the product formula was extended to plates bounded by the graph of a positive function. Develop the continuous analogue of Theorems~\ref{thm:steiner} and~\ref{thm:twosided}, and compare the convergence of the centres to the right end point with the method of Laplace.
\end{problem}

\section{Remarks on earlier papers}\label{sec:earlier}

We correct some statements in the earlier papers \cite{Rahimi09,Rahimi17,MehdiNezhadRahimi}.
\begin{enumerate}[\rm(a)]
\item In \cite[Remark~1.1 and Example~1.3]{MehdiNezhadRahimi} the position of the centre was inferred from the number of points in each half of the interval $[x_1,x_n]$. This is not correct: three of the five points $1,49,50,99,100$ lie to the left of the midpoint $50.5$, but their centre is $59.8$. The statements of \cite[Theorem~2.11 and Proposition~2.13]{MehdiNezhadRahimi} remain true, and follow from Lemma~\ref{lem:lever}. Lemma~2.4 of \cite{MehdiNezhadRahimi} requires $n\ge4$, since $2f_1=f_3$. Also, the bounds of \cite{MehdiNezhadRahimi} are weaker than the elementary estimate $N<2f_{n_k}$, valid for every Fibonacci representation $N=f_{n_1}+\dots+f_{n_k}$ in the sense of that paper, since the terms other than $f_{n_k}$ are distinct Fibonacci numbers not exceeding $f_{n_k-2}$, and their sum is at most $f_{n_k}-2$.
\item In the first assertion of \cite[Theorem~2.9(a)]{Rahimi17} the final factor $y$ must be replaced by $y^{m/r}$; as stated, the assertion fails for $1^3+2^2=5$. Corollary~2.11 of \cite{Rahimi17} applies Beal's conjecture to equations with coefficients, to which it does not apply; for instance $19\cdot1^3+2^3=3^3$. The proof of \cite[Lemma~2.5]{Rahimi17} contains a misprint ($+$ for $=$); the statement is correct, and is contained in Theorem~\ref{thm:steiner}.
\item In \cite{Rahimi09} there are only misprints; Theorem~3 there is correct.
\end{enumerate}

\end{document}